\documentclass[number,secthm,dvips]{arxbj}
\usepackage{graphicx}

\aid{0}
\volume{14}
\issue{1}
\pubyear{2008}
\firstpage{155}
\lastpage{179}
\doi{10.3150/07-BEJ6150}

\makeatletter

\newtheorem{lm}[thm]{Lemma}

\newtheorem{proposition}[thm]{Proposition}

\theoremstyle{definition}

\newremark{example}[thm]{Example}
\newremark{remark}[thm]{Remark}
\makeatother

\begin{document}
\begin{frontmatter}

\title{Sequential Monte Carlo
smoothing with application to parameter estimation in nonlinear
state space models}
\runtitle{Sequential Monte Carlo smoothing}

\begin{aug}
\author[1]{\fnms{Jimmy} \snm{Olsson}\corref{}\thanksref{1,e1}\ead[label=e1,mark]{olsson@enst.fr}},
\author[1]{\fnms{Olivier} \snm{Capp\'e}\thanksref{1,e2}\ead[label=e2,mark]{cappe@enst.fr}},
\author[2]{\fnms{Randal} \snm{Douc}\thanksref{2}\ead[label=e3]{douc@cmapx.polytechnique.fr}} and
\author[1]{\fnms{\'Eric}~\snm{Moulines}\thanksref{1,e4}\ead[label=e4,mark]{moulines@enst.fr}}
\runauthor{J. Olsson et al.}
\pdfauthor{Jimmy Olsson, Olivier Cappe, Randal Douc and Eric Moulines}
\address[1]{Ecole Nationale Sup\'{e}rieure des
T\'{e}lecommunications, France. \printead{e1};\break \printead
*{e2}; \printead*{e4}}
\address[2]{CMAP, Ecole Polytechnique, France. \printead{e3}}
\end{aug}

% HISTORY:
\received{\smonth{9} \syear{2006}}
\revised{\smonth{9} \syear{2007}}

% ABSTRACT
%
\begin{abstract}
This paper concerns the use of sequential Monte Carlo methods (SMC) for
smoothing in general state space models. A well-known problem when
applying the standard SMC technique in the smoothing mode is that the
resampling mechanism introduces degeneracy of the approximation in the
path space. However, when performing maximum likelihood estimation via
the EM
algorithm, all functionals involved are of additive form for a
large subclass of models. To cope with the problem in this case, a
modification of the standard method (based on a technique proposed by
Kitagawa and Sato) is suggested. Our
algorithm relies on forgetting properties of the filtering dynamics
and the quality of the estimates produced is investigated, both
theoretically and via simulations.
\end{abstract}

% KEYWORDS
%

\end{frontmatter}
%
%s1 ###
\section{Introduction}\label{s1}

In this paper, we study SMC methods for smoothing in nonlinear state
space models. We consider a bivariate process $(X, Y)$, where $X
\triangleq\{X_k ; k \geq0 \}$ is a homogeneous discrete-time Markov
chain taking values in some state space $(\mathsf{X}, \mathcal{X})$. We
let $(Q_{\theta}, \theta\in\Theta\subseteq\mathbb{R}^d)$ and $\nu$
denote the Markov transition kernel and the initial distribution of $X$,
respectively. The family $\{Q_{\theta}(x,\cdot)$; $x \in
\mathsf{X}$, $\theta\in\Theta\}$ is assumed to be dominated by the
probability
measure $\mu$ on $(\mathsf{X}, \mathcal{X})$ and we denote by
$q_{\theta}
(x,\cdot)$ the corresponding Radon--Nikodym derivatives. In this
framework, $X$ is not observed and measurements must be made
through the process $Y \triangleq\{ Y_k ; k \geq0 \}$ taking values in
some measurable space $(\mathsf{Y}, \mathcal{Y})$. These observed
variables are
conditionally independent, given the sequence $\{X_k ; k \geq0 \}$,
and the conditional distribution of $Y_k$ depends
only on $X_k$. We denote by $\mathcal{G}_{k}$ the
$\sigma$-algebra generated by the observed process from time zero to
$k$. Furthermore, there exist, for all $x \in\mathsf{X}$ and
$\theta\in\Theta$, a density function $y \mapsto g_{\theta}(x,y)$ and
a measure $\lambda$ on $(\mathsf{Y}, \mathcal{Y})$ such that, for $k
\geq0$,
\[
\operatorname{\mathbb{P}}_\theta( Y_k \in A
| X_k = x) = \int_A g_{\theta}(x, y)
\lambda(\mathrm{d}y)  \qquad\mbox{for all } A \in\mathcal{Y}.
\]
Here, we have written $\operatorname{\mathbb{P}}_\theta$ for the law
of the bivariate
Markov chain $\{(X_k, Y_k); k \geq0\}$ under the model
parameterized by $\theta\in\Theta$ and we denote by $\mathbb
{E}_\theta$
the associated expectation.

For $i \leq j$, let $X_{i:j} \triangleq(X_i, \ldots, X_j)$;
similar vector notation will be used for other quantities. In many
situations, it is required to compute expectation values of the form
$\mathbb{E}_{\theta}[t_n(X_{0:n})| \mathcal{G}_{n} ]$, where $t_n$
is a
real-valued, measurable function. In this paper, we focus on the case
where $t_n$ is an \textit{additive functional} given by
%
%e1 ###
%
\begin{equation} \label{eqtndef}
t_n(x_{0:n}) = \sum_{k=0}^{n-1} s_k(x_{k:k+1}),
\end{equation}
where $\{s_k; k \geq0 \}$ is a sequence of measurable functions
(which may depend on the observed values $Y_{0:n}$).

As an example of when smoothing of such additive functionals is
important, consider the case of maximum likelihood estimation via the
\textit{EM algorithm}. Having an initial estimate $\theta'$ of the parameter
vector available, an improved estimate is obtained (we refer
to Capp\'e \textit{et} \textit{al.} \cite{cappemoulinesryden2005},
Section~10.2.3)
by means of computation and maximization of $\mathcal{Q}(\theta;
\theta')$
with respect to $\theta$, where $\mathcal{Q}(\theta; \theta')$ is
defined by
\begin{eqnarray*}
\mathcal{Q}(\theta; \theta') &\triangleq&\mathbb{E}_{\theta'} \Biggl[
 \sum_{k=0}^{n-1} \log q_{\theta} (X_k, X_{k+1}) \Big|
\mathcal{G}_{n} \Biggr] + \mathbb{E}_{\theta'}
\Biggl[ \sum_{k=0}^n \log g_{\theta} (X_k, Y_k) \Big|
\mathcal{G}_{n}
\Biggr] \\ &&{}+ \mathbb{E}_{\theta'}[  \log\nu
(X_0) | \mathcal{G}_{n}].
\end{eqnarray*}
This procedure is recursively repeated in order to obtain
convergence to a stationary point $\theta_\star$ of the
\textit{log-likelihood function} $\ell_{\nu, n} (\theta;
Y_{0:n}) \triangleq\log\mathrm{L}_{\nu, n}(\theta;
Y_{0:n})$, where, for $y_{0:n} \in\mathsf{Y}^{n+1}$,
\[
\mathrm{L}_{\nu, n}(\theta; y_{0:n}) \\ \triangleq
\int_{\mathsf{X}^{n+1}} g_{\theta} (x_0,y_0) \nu(x_0)
\prod_{k=1}^n q_{\theta} (x_{k-1}, x_k) g_{\theta} (x_k, y_k)
\mu^{\otimes(n+1)}(\mathrm{d}x_{0:n}).
\]
The computation of smoothed sum functionals of the above form will
also be the key issue when considering direct maximum likelihood
estimation via the \textit{score function} $\nabla_\theta\ell_{\nu
, n}(\theta; y_{0:n})$; again see Capp\'{e} \textit{et} \textit{al.}
(\cite{cappemoulinesryden2005}, Section~10.2.3)
for details.

By applying Bayes' formula, it is straightforward to derive
recursive formulas for expectations of the additive type discussed
above. However, tractable closed form solutions are available only if
the state space $\mathsf{X}$ is finite or the model is linear and Gaussian.

SMC methods (also known as \textit{particle filtering methods}) constitute
a class of algorithms that are well suited for providing approximate
solutions of the smoothing and filtering recursions. In recent years, SMC
methods have been applied, sometimes very successfully,
in many different fields (see Doucet \textit{et al.} \cite{doucetdefreitasgordon2001}
and Ristic \textit{et al.} \cite{risticarulampalamgordon2004} and the references
therein). A well-known problem when applying SMC methods to
sample the joint smoothing distribution is that the resampling
mechanism of the
particle filter introduces degeneracy of the particle
trajectories. Doucet \textit{et al.} \cite{doucetgodsillwest2004}
suggest a procedure where this is avoided through an additional
resampling pass in the time-reversed direction. The
resulting algorithm is well suited to sample from the joint smoothing
distribution, but appears unnecessarily complex, computationally, for
approximating additive smoothing functionals of the form~(\ref{eqtndef}).\looseness=1

In this paper, we study an SMC technique to smooth additive
functionals based on a fixed-lag smoother presented by
Kitagawa and Sato \cite{kitagawasato2001}. The method exploits the \textit{forgetting
properties} on the conditional hidden chain and is
not affected by the degeneracy of the particle trajectories. Compared
to  Doucet \textit{et al.} \cite{doucetgodsillwest2004}, computational requirements are
marginal. Furthermore, we perform, under suitable regularity
assumptions on the latent chain, a theoretical analysis of the
behavior of the estimates obtained. It turns out that the $\mathsf{L}^{p}$
error and bias are upper bounded by  quantities proportional to
$n \log n / \sqrt{N}$ and $n \log n / N$, respectively, where $N$
denotes the
number of particles and $n$ the number of observations.

In  comparison, applying the results of
Del~Moral and Doucet (\cite{delmoraldoucet2003}, Theorem~4) to a
functional of type (\ref{eqtndef}) provides
a bound proportional to $n^2 / \sqrt{N}$
on the $\mathsf{L}^{p}$ error for the standard trajectory-based particle
smoother. Finally, we apply, for a noisily observed autoregressive
model and the stochastic volatility model proposed by
Hull and White \cite{hullwhite1987}, the technique to the \textit{Monte Carlo EM} (MCEM)
\textit{algorithm} (Wei and Tanner \cite{weitanner1991}).

%s2 ###
\section{Particle approximation of additive functionals}
\label{sectionParticleApproximation}
%s2.1 ###
\subsection{The smoothing recursion} The \textit{joint smoothing
distribution} $\phi_{\nu, 0:n|n}$ is the probability measure defined,
for $A \in\mathcal{X}^{\otimes(n+1)}$, by
\[
\phi_{\nu, 0:n|n}[Y_{0:n}](A;\theta) \triangleq
\operatorname{\mathbb{P}}_\theta( X_{0:n} \in A
| \mathcal{G}_{n}).
\]
Under the assumptions above, the joint smoothing distribution has a
density (for which we will use the same symbol) with respect to
$\mu^{\otimes(n+1)}$ satisfying, for all $y_{0:k+1} \in
\mathsf{Y}^{k+2}$, the recursion
%
%e2 ###
%e1 ###
%
\begin{eqnarray} \label{smoothingrecursion}
&&\phi_{\nu, 0:k+1|k+1}[y_{0:k+1}](x_{0:k+1};
\theta)\nonumber\\
&&\quad= \frac{\mathrm{L}_{\nu, k} (\theta;
y_{0:k})}{\mathrm{L}_{\nu, k+1}(\theta; y_{0:k+1})}
q_{\theta} (x_k, x_{k+1}) g_{\theta} (x_{k+1},y_{k+1})
\phi_{\nu, 0:k|k}[y_{0:k}](x_{0:k}; \theta).
\end{eqnarray}
For notational conciseness, we will omit the explicit
dependence on the observations from the notation for the smoothing
measure and replace
$\phi_{\nu, 0:k|k}[y_{0:k}](\cdot ;\theta)$ by
$\phi_{\nu, 0:k|k}(\cdot ;\theta)$.

Particle filtering, in its most basic form, consists of approximating
the exact smoothing relations by propagating particle trajectories in
the state space of the hidden chain. Given a fixed sequence of
observations, this is done according to the following scheme. In order to
keep the notation simple, we fix the model parameters and omit
$\theta$ from the notation throughout this part.

At time zero, $N$ random variables $\{\xi_{0}^{N,i} ; 1 \leq i \leq N
\}$ are drawn from a common probability measure $\varsigma$ such that
$\nu\ll\varsigma$. These \textit{initial particles} are assigned
the \textit{importance weights}
$\omega_{0}^{N, i} \triangleq W_0 (\xi_{0}^{N,i})$, $1 \leq i
\leq N$, where, for $x \in\mathsf{X}$, $W_0(x) \triangleq g_{}(x,y_0)
\, \mathrm{d}
\nu/\mathrm{d}\varsigma(x)$, providing $\sum_{i=1}^N \omega
_{0}^{N, i}
f (\xi_{0}^{N,i})/ \sum_{i=1}^N \omega_{0}^{N, i}$ as an
importance sampling estimate of $\phi_{\nu, 0|0}f$ for $f \in
\mathcal{B}_{\mathrm{b}}(\mathsf{X}^{})$. Henceforth, the particle
paths $\xi_{0:m}^{N,i} \triangleq
[\xi_{0:m}^{N,i}(0), \ldots, \xi_{0:m}^{N,i}(m)]$, $1 \leq i \leq
N$, are
recursively updated according to the following procedure.

At time $k$, let $\{ (\xi_{0:k}^{N,i},\omega_{k}^{N, i}); 1 \leq i
\leq N \}$ be
a set of
weighted particles approximating $\phi_{\nu, 0:k|k}$, in the sense that
$\sum_{i=1}^N \omega_{k}^{N, i} f(\xi_{0:k}^{N,i}) / \Omega_{k}^N$, with
$\Omega_{k}^N \triangleq\sum_{i=1}^N \omega_{k}^{N, i}$ and $f \in
\mathcal{B}_{\mathrm{b}}(\mathsf{X}^{k+1})$, is an estimate of the
expectation
$\phi_{\nu, 0:k|k}f$. Then, an updated weighted
sample $\{(\xi_{0:k+1}^{N,i},\omega_{k+1}^{N, i}); 1 \leq i \leq N \}$,
approximating the distribution $\phi_{\nu, 0:k+1|k+1}$, is
obtained by, first, simulating $\xi_{0:k+1}^{N,i} \sim
R_{k}^\mathrm{p}(\xi_{0:k}^{N,i},\cdot)$, where the kernel
$R_{k}^\mathrm{p}$ is of type $R_{k}^\mathrm{p}(x_{0:k},f) =
\int_\mathsf{X}f(x_{0:k},x_{k+1})  R_k(x_k, \mathrm{d}x_{k+1})$, with
$f \in\mathcal{B}_{\mathrm{b}}(\mathsf{X}^{k+2})$ and each $R_k$
being a Markov transition
kernel. The new particles are simulated independently of each other
and the
special form of $R_{k}^\mathrm{p}$ implies that past particle
trajectories are kept unchanged throughout this \textit{mutation
step}. A popular choice is to set $R_{k} \equiv Q$, yielding the
so-called \textit{bootstrap filter}; more sophisticated techniques
involve proposals depending on the observed values (see
Example \ref{examplestovol2}). Second, when the observation
$Y_{k+1} = y_{k+1}$ is available, the importance weights are updated
according to $\omega_{k+1}^{N, i} = \omega_{k}^{N, i}
W_{k+1}[\xi_{0:k+1}^{N,i}(k:k+1)]$, where, for $(x,x') \in
\mathsf{X}^2$, $ W_k(x,x') \triangleq
g_{}(x', y_k)  \,\mathrm{d}Q(x, \cdot) / \mathrm{d}R_{k-1}(x,
\cdot)(x')$. Now, for $f \in\mathcal{B}_{\mathrm{b}}(\mathsf
{X}^{k+2})$, the self-normalized
estimate $\phi_{\nu, 0:k+1|k+1}^Nf \triangleq\sum_{j=1}^N
\omega_{k+1}^{N, j} f(\xi_{0:k+1}^{N,j}) / \Omega_{k+1}^N$ provides an
approximation of $\phi_{\nu, 0:k+1|k+1}$.

To prevent degeneracy, a \textit{resampling mechanism} is introduced. In
its simpler form, resampling amounts to drawing, conditionally
independently, indices $I_k^{N,1}, \dots, I_k^{N,N}$ from the set $\{
1,\dots,N\}$,
multinomially with respect to the normalized weights $\omega
_{k}^{N, j} /
\Omega^{N}_{k}$, $1 \leq j \leq N$. Now, a~new
equally weighted sample $\{
\hat\xi_{0:k}^{N,i}; 1 \leq i \leq N \}$ is
constructed by setting $ \hat{\xi}_{0:k}^{N,j} = \xi_{0:k}^{N,I_{k}^{N,j}}$.
After the resampling procedure, the weights are all reset as
$\omega_{k}^{N, i}= 1/N$, yielding another estimate, $\widehat{\phi
}_{\nu, 0:k|k}^Nf
\triangleq\sum_{i=1}^N f(\hat\xi_{0:k}^{N,i}) / N$, of
$\phi_{\nu, 0:k|k}$. Note that the resampling mechanism might modify
the whole trajectory of a certain particle, implying that, in general,
for $m \leq n$, $\xi_{0:n}^{N,i}(m) \neq\xi_{0:n+1}^{N,i}(m)$. The
multinomial resampling method is not the only conceivable way to carry
out the selection step (see e.g.   Doucet \textit{et al.} \cite{doucetdefreitasgordon2001}).

Using the weighted samples $\{ (\xi_{0:k}^{N,j}, \omega_{k}^{N, j}); 1
\leq j \leq N \}$, $0 \leq k \leq n$, produced under the parameter
$\theta\in\Theta$, an approximation of $\gamma_{\theta, n}
\triangleq
\mathbb{E}_\theta[t_n(X_{0:n}) | \mathcal{G}_{n}]$ is obtained by constructing
the estimators
%
%e3 ###
%
\begin{equation} \label{eqsmoothingest}
\gamma_{\theta, n}^N = \frac{1}{\Omega_{n}^N} \sum_{j=1}^N
\omega_{n}^{N, j} t_n(\xi_{0:n}^{N,j}) \quad\mbox{or} \quad
\widehat{\gamma}_{\theta, n}^N = \frac{1}{N} \sum_{j=1}^N t_n (
\hat{\xi}_{0:n}^{N,j}).
\end{equation}

When the functional $\{t_n\}$ has the form given in (\ref{eqtndef}),
it is
straightforward to verify that recording all of the particle trajectories is
indeed not required to evaluate (\ref{eqsmoothingest}): upon defining
$t_k^{N,i}
\triangleq t_k ( {\xi}_{0:k}^{N,j})$, we have, for $k\geq1$,
%
%e4 ###
%
\begin{equation} \label{equpdateapprox}
t_{k+1}^{N,i} =
\cases{t_k^{N,i} + s_k
[{\xi}_{0:k+1}^{N,i}(k:k+1) ], \quad & if no
resampling occurs,\cr
t_k^{N,I_{k+1}^{i}} + s_k
[\hat{\xi}_{0:k+1}^{N,i}(k:k+1)], \quad & if
resampling occurs. }
\end{equation}

The recursion is initialized by $t^{N,i}_1 = t_1(\xi_{0:1}^{N,i})$.
In accordance with (\ref{eqsmoothingest}), $\gamma_n^N$ is
obtained as $\sum_{i=1}^N \omega_{n}^{N, i}  t_n^{N,i} / \Omega
_{n}^N$. Hence, for
each particle ${\xi}_{0:k}^{N,i}$, we need only record its current
position ${\xi}_{0:k}^{N,i}(k)$, weight $\omega_{k}^{N, i}$ and associated
functional
value $t_k^{N,i}$. Thus, the method necessitates only minor
adaptations once the particle filter has been implemented.

\begin{figure}%

\includegraphics{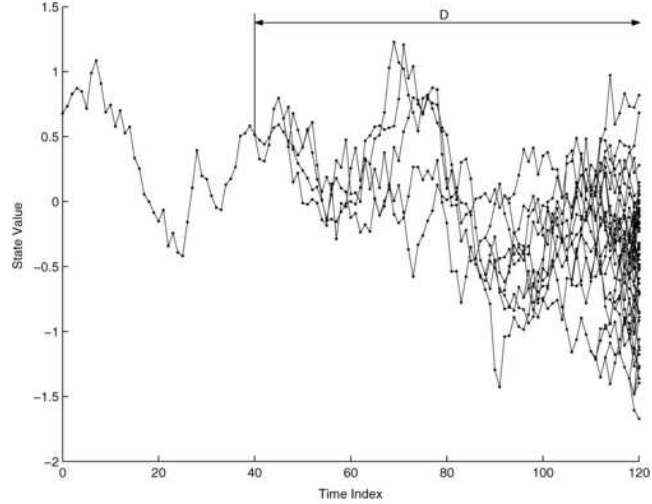}

\caption{Typical particle trajectories for $N=50$; see Section
\protect\ref{sectionApplicationsMLE} for details regarding model and algorithm.}
\label{figpath}
\end{figure}

As illustrated in Figure~\ref{figpath}, as $n$ increases, the path
trajectories system collapses,
and the estimators (\ref{eqsmoothingest}) are not reliable for
sensible $N$ values (see Doucet \textit{et al.} \cite{doucetdefreitasgordon2001},
Kitagawa and Sato \cite{kitagawasato2001} and Andrieu and Doucet \cite{andrieudoucet2003} for a
discussion).

%f1 ###
%

To cope with this drawback, we suggest the following method,
based on a
technique proposed by Kitagawa and Sato \cite{kitagawasato2001}. By the
forgetting property of the time-reversed conditional hidden chain
(Theorem \ref{mixingconditionalchainsth}), we
expect that, for a large enough integer $\Delta_n \leq n-k$,
%e5 ###
%
\begin{equation}\label{keyforgettingrelation}
\mathbb{E}_\theta[  s_k (X_{k:k+1}) |
\mathcal{G}_{n} ]
\approx\mathbb{E}_\theta[  s_k (X_{k:k+1}) |
\mathcal{G}_{k+\Delta_n} ],
\end{equation}
yielding, with $k(\Delta_n) \triangleq(k+\Delta_n) \wedge n$,
\[
\gamma_{\theta, n} = \mathbb{E}_\theta\Biggl[  \sum
_{k=0}^{n-1} s_k
(X_{k:k+1}) \Big|
\mathcal{G}_{n} \Biggr]
\approx\sum_{k=0}^{n-1} \mathbb{E}_\theta\bigl[  s_k
(X_{k:k+1}) |
\mathcal{G}_{k(\Delta_n)} \bigr].
\]
The above relation suggests that waiting for all of the trajectories to
collapse -- as (\ref{equpdateapprox}) implies -- is not
convenient. Instead, when the particle population
$N$ is sufficiently large so that (\ref{keyforgettingrelation}) is
valid for a lag $\Delta_n$ which may be far smaller than the typical
collapsing time, one should apply the two approximations
\begin{eqnarray}
\label{eqfixedlagsmoothingestimateweighted}
\gamma_{\theta, n}^{N,\Delta_n} &\triangleq&\sum_{k=0}^{n-1}
\sum_{j=1}^N \frac{\omega_{k(\Delta_n)}^{N, j}}{\Omega_{k(\Delta
_n) }^N} s_k
\bigl[\xi_{0:k(\Delta_n)}^{N,j}(k\dvtx k+1) \bigr], \\
\label{eqfixedlagsmoothingestimateresamp}
\widehat{\gamma}_{\theta, n}^{N,\Delta_n} &\triangleq&\frac{1}{N}
\sum_{k=0}^{n-1}
\sum_{j=1}^N s_k \bigl[\xi_{0:k(\Delta_n)}^{N,j}(k\dvtx k+1) \bigr]
\end{eqnarray}
of $\gamma_{\theta, n}$. Although somewhat
more involved than the standard approximation (\ref{eqsmoothingest}),
the above lag-based approximation  may be updated recursively by
recording the recent history of the particles as well as
the accumulated contribution of terms that will no longer get updated.
Thus, apart from increased storage requirements, computing the
lag-based approximation $\widehat{\gamma}^{N,\Delta_n}_{\theta, n}$
is clearly
not, from a computational point of view, more demanding than computing~$\widehat{\gamma}_{\theta, n}^N$.\looseness=1

%s3 ###
\section{Theoretical evaluation of the fixed-lag technique}
\label{sectiontheoretics}
\setcounter{equation}{0}
To accomplish the robustification above, we
need to specify the lag $\Delta_n$ and how this lag should depend on
$n$. This is done by examining the quality
of the estimates produced by the algorithm in terms of bias and
$\mathsf{L}^{p}$
error. Of particular interest is how these errors are affected
by the lag and whether it makes their dependence on $n$ and $N$ more
favorable in comparison with the standard trajectory-based approach.%\looseness=1

The validity of \ref{keyforgettingrelation} is based on the
assumption that the conditional hidden chains -- in the forward as well as
the backward directions -- have \textit{forgetting properties},
that is, the distributions of two versions of each chain starting at
different initial distributions approach each other as time increases. This
property depends on the following uniform ergodicity
conditions on the model, which imply that forgetting occurs at a
\textit{geometrical} rate:%\looseness=1

\begin{enumerate}[(A1)]
\item[(A1)]

\begin{itemize}[(i)]
\item[(i)] $\sigma_-\triangleq\inf_{\theta\in\Theta} \inf
_{x,x' \in
\mathsf{X}}
q_{\theta} (x,x')>0$, $\sigma_+\triangleq\sup_{\theta\in
\Theta} \sup_{x,x' \in
\mathsf{X}} q_{\theta} (x,x') < \infty$;
\item[(ii)] for all $y \in\mathsf{Y}$, $\sup_{\theta\in
\Theta}
\| g_{\theta} (\cdot, y) \|_{\mathsf{X}^{}, \infty} <
\infty$,
$\inf_{\theta\in\Theta} \int_\mathsf{X}g_{\theta} (x,y)  \mu(
\mathrm{d}
x) > 0$.
\end{itemize}
\end{enumerate}
Under (A1), we define
%
%e1 ###
%
\begin{equation} \label{eqrhodef}
\rho\triangleq1 - \frac{\sigma_-}{\sigma_+}.
\end{equation}

We now define the Markov transition kernels that generate the conditional
hidden chains. For any two transition kernels $K$ and $T$ from
$(\mathsf{E}_1, \mathcal{E}_1)$ to $(\mathsf{E}_2, \mathcal{E}_2)$ and
$(\mathsf{E}_2, \mathcal{E}_2)$ to $(\mathsf{E}_3, \mathcal{E}_3)$,
respectively, we define the product transition kernel by $KT(x, A)
\triangleq\int_{\mathsf{E}_2} T(z, A) K(x, \mathrm{d}z)$ for $x
\in
\mathsf{E}_1$ and $A \in\mathcal{E}_3$.

Introduce, for $f \in\mathcal{B}_{\mathrm{b}}(\mathsf{X}^{k+2})$,
$x_{0:k} \in\mathsf{X}^{k+1}$
and $y_{k+1} \in\mathsf{Y}$, the unnormalized pathwise transition
kernel $L_{k} (x_{0:k}, f ; \theta) \triangleq\int_\mathsf{X}
f(x_{0:k+1}) g_{\theta}(x_{k+1}, y_{k+1})  Q_{\theta}
(x_k, \mathrm{d}x_{k+1})$. Assumption~(A1)
makes this
integral well defined for all $k \geq0$. We will often consider compositions
\[
L_{k} \cdots L_{m} (x_{0:k}, f ; \theta) = \int_{\mathsf{X}^{m-k+1}}
f(x_{0:m+1}) \prod_{i=k}^m [ g_{\theta}(x_{i+1}, y_{i+1})
 Q_{\theta}(x_i, \mathrm{d}x_{i+1}) ]
\]
with $f \in\mathcal{B}_{\mathrm{b}}(\mathsf{X}^{m+2})$, $x_{0:k}
\in\mathsf{X}^{k+1}$ and
$y_{0:k} \in\mathsf{Y}^{m-k+1}$, and it is clear that, for all
$k \leq m$, the function $L_{k} \cdots L_{m} (x_{0:k},
\mathsf{X}^{m+2} ; \theta)$ depends only on $x_k$. Thus, a version of this
function comprising only the last component is well defined and we
write $L_{k} \cdots L_{m} (x_k, \mathsf{X}^{m+2} ; \theta)$ in this case.
For $k > m$, we set $L_{k} \cdots L_{m} \equiv
\operatorname{Id}$. Using this notation and given $n \geq0$, the
\textit{forward smoothing kernels} given by, for $k \geq0$, $x_k \in
\mathsf{X}$ and $A \in\mathcal{X}$, $\mathrm{F}_{k|n} (x_k, A ;
\theta) \triangleq
\operatorname{\mathbb{P}}_\theta(  X_{k+1} \in A
| X_k = x_k,
\mathcal{G}_{n})$, can, for indices $0 \leq k < n$ and $y_{k+1}
\in\mathsf{Y}$, be written as
%
%e4 ###
%e3 ###
%e2 ###
%
\begin{eqnarray} \label{eqforwardkernel}
&&\mathrm{F}_{k|n} (x_k, A ; \theta) \nonumber\\[-8pt]\\[-8pt]
\nonumber&&\quad=
\cases{ \displaystyle\int_A
\dfrac{g_{\theta}(x_{k+1}, y_{k+1})
L_{k+1} \cdots L_{n-1} (x_{k+1}, \mathsf{X}^{n+1}; \theta)
Q_{\theta}(x_k, \mathrm{d}x_{k+1})}{L_{k} \cdots L_{n-1} (x_k, \mathsf
{X}^{n+1}; \theta)}, \cr
\hspace*{47pt}\qquad\mbox{ for }\     0 \leq k < n, \cr
Q_{\theta}(x_k, A),
\qquad\mbox{ for  }\ k \geq n.}
\end{eqnarray}
Analogously, for the time-reversed conditional hidden chain, we consider the
\textit{backward smoothing kernels} defined by
$\mathrm{B}_{\nu, k|n} (x_{k+1}, A ; \theta) \triangleq
\operatorname{\mathbb{P}}_\theta
( X_k \in A | X_{k+1} = x_{k+1}, \mathcal{G}_{n}
)$, where $k \geq0$, $x_{k+1} \in\mathsf{X}$ and $A \in
\mathcal{X}$.
Note that
$\mathrm{B}_{\nu, k|n}$ depends on the initial distribution of the
latent chain. The
backward kernel can be expressed as
\begin{eqnarray*}
&&\mathrm{B}_{\nu, k|n} (x_{k+1}, A ; \theta)
\\&&\quad =
\cases{\dfrac{
\int_A q_{\theta} (x_k,
x_{k+1})  \phi_{\nu,k} (\mathrm{d}x_k; \theta)}{\int_\mathsf{X}
q_{\theta}(x'_k, x_{k+1})  \phi_{\nu,k}(\mathrm{d}x'_k;
\theta)}, &\quad for  $0 \leq k \leq n$, \cr
\dfrac{ \int_A \int_\mathsf{X}
q_{\theta}(x_k, x_{k+1}) q_{\theta}^{k-n}(x_n, x_k)
\phi_{\nu,n}(\mathrm{d}x_n; \theta)  \mu(\mathrm{d}
x_k)}{\int_\mathsf{X}q_{\theta}^{k-n+1}(x'_n,
x_{k+1})  \phi_{\nu,n}(\mathrm{d}x'_n; \theta)}
, &\quad for  $k > n$,}
\end{eqnarray*}
where, for $m \geq1$, $q_{\theta}^m$ denotes the density of the
$m$-step kernel $Q_{\theta}^m$.

The following theorem (see Del~Moral \cite{moral2004}, page~143), stating
geometrical ergodicity of the forward and backward chains, is
instrumental for the developments which are to follow.
\begin{thm}
\label{mixingconditionalchainsth}
Assume \textup{(A1)} and let $\rho$ be defined in
(\ref{eqrhodef}). Then, for all $k \geq m \geq0$, all $\theta\in
\Theta$, all probability measures $\nu_1$, $\nu_2$ on $\mathcal
{X}$ and
all $y_{0:n} \in\mathsf{Y}^{n+1}$,
\begin{eqnarray*}
\| \nu_1 \mathrm{F}_{m|n} \cdots\mathrm{F}_{k|n} (\cdot ;
\theta) - \nu_2
\mathrm{F}_{m|n} \cdots\mathrm{F}_{k|n} (\cdot; \theta) \|
_{\mathrm{TV}}
&\leq&\rho^{k-m+1}, \\
\| \nu_1 \mathrm{B}_{\nu, k|n} \cdots\mathrm{B}_{\nu, m|n}
(\cdot ; \theta) - \nu_2
\mathrm{B}_{\nu, k|n} \cdots\mathrm{B}_{\nu, m|n} (\cdot;
\theta) \|_{\mathrm{TV}}
&\leq&\rho^{k-m+1}.
\end{eqnarray*}
\end{thm}

Assumption~(A1) typically requires that
$\mathsf{X}$ is a
compact set, but some very recent papers
(Douc~et al. \cite{doucfortmoulinespriouret2007}, Kleptsyna and Veretennikov  \cite{kleptsynaveretennikov2007})
provide results that establish geometric forgetting under considerably
weaker assumptions. Applying these results within our framework would,
however, make the analysis far more complicated since the
provided bounds are uniform neither in the observations nor the initial
distributions.

%s3.1 ###
\subsection{Main results}
For the sake of simplicity, let us assume that multinomial resampling is
applied at every iteration. Moreover, let the observations used by the particle
filter be generated by a state space model with kernel, measurement
density and initial distribution $\bar{Q}$, $\bar{g}$ and
$\bar{\nu}$, respectively. We stress that $\bar{Q}$
and $\bar{g}$ are not assumed to belong to the parametric family
$\{(Q_\theta,g_\theta); \theta\in\Theta\}$. Using these observed
values as input, the evolution of the particle cloud follows the usual
dynamics ($Q_\theta$,$g_\theta$,$\nu$, $\theta\in\Theta$)
and, in this setting, it is easily verified that the process
$\{Z_k; k \geq0\}$, with $Z_k \triangleq[\xi_{0:k}^{N,1}(k-1\dvtx k), \ldots,
{\xi_{0:k}^{N,N}}(k-1\dvtx k), X_k, Y_k]$, is a Markov chain on
$\mathsf{X}^{2 N +1} \times\mathsf{Y}$. We denote by
$\bar{\mathbb{P}}^N_\theta$
and
$\bar{\mathbb{E}}^N_\theta$
the law of this
chain and the associated expectation, respectively, and define the
filtration $\{ \mathcal{F}_{k}^N ; k \geq0 \}$ by $\mathcal{F}_{k+1}^N
\triangleq\mathcal{F}_{k}^N \vee\sigma({\xi_{0:k+1}^{N,1}}, \ldots,
{\xi_{0:k+1}^{N,N}})$
with $\mathcal{F}_{0}^N \triangleq\sigma(
\xi_{0}^{N,1}, \ldots, \xi_{0}^{N,N})$. The marginal of
$\bar{\mathbb{P}}^N_\theta$
with respect to
$\{ (X_k, Y_k) ; k \geq0 \}
$ and the
associated expectation are denoted by
$\bar{\mathbb{P}}$ and~$\bar{\mathbb{E}}$, respectively. For any
integer $p \geq
1$, random variable $V \in\mathsf{L}^{p}(\bar{\mathbb{P}}^N_\theta)$ and
sub-$\sigma$-algebra $\mathcal{A} \subseteq\sigma(\{Z_k; k \geq0\})$
we define the conditional $\mathsf{L}^{p}$ norm $\| V \|_{p | \mathcal{A}}
\triangleq
(\bar{\mathbb{E}}^N_\theta [|V|^p | \mathcal{A} ])^{1/p}$.
\begin{enumerate}
\item[(A2)]For all $k \geq0$, $\theta\in\Theta$ and $y_k \in\mathsf{Y}$,
$\| W_k(\cdot; \theta) \|_{\mathsf{X}^{2}, \infty} <
\infty$.\end{enumerate}
\begin{remark}
In case of the bootstrap particle filter, for which $R_{\theta, k}
\equiv Q_\theta$, assumption \textup{(A2)} is
implied by
assumption~\textup{(A1)}. The same is true for the so-called
\textit{optimal kernel} used in Example~\ref{examplestovol2}.
\end{remark}
\begin{thm} \label{mainresulttheorem}
Assume \textup{(A1)} and
\textup{(A2)}. There then exist universal
constants $B_p$ and $B$, $B_p$ depending only on $p$, such that the
following holds true for all $n \geq0$, $\theta\in\Theta$, $\Delta
_n \geq
0$ and $N \geq1$:
\begin{itemize}
\item[(i)] for all $p \geq2$,
\begin{eqnarray*}
&&\| \widehat{\gamma}^{N,\Delta_n}_{\theta,n} - \gamma_{\theta,n}
\|_{p| \mathcal{G}_{n}}\\
&&\quad \leq2 \rho^{\Delta_n}
\sum_{k=0}^{n-\Delta_n}
\| s_k \|_{\mathsf{X}^{2}, \infty} \\
&&\quad\quad{}+
\frac{B_p}{\sqrt{N}(1-\rho)} \sum_{k=0}^{n-1} \| s_k \|
_{\mathsf{X}^{2}, \infty} \Biggl[
\frac{1}{\sigma_-} \sum_{m=1}^{k(\Delta_n)}
\frac{\| W_m(\cdot; \theta) \|_{\mathsf{X}^2, \infty} \rho^{0
\vee
(k-m)}}{\mu g_{\theta}(Y_m)}\\
&&\hspace*{183pt}{}+
\frac{\| W_0(\cdot; \theta) \|_{\mathsf{X}, \infty}
\rho^k}{\nu g_{\theta}(Y_0)} + 1 \Biggr];
\end{eqnarray*}
\item[(ii)]
\begin{eqnarray*}
&&\bigl| \bar{\mathbb{E}}^N_\theta [
 \widehat{\gamma}^{N,\Delta_n}_{\theta,n} |
\mathcal{G}_{n} ] - \gamma_{\theta,n} \bigr|\\
 &&\quad\leq2 \rho
^{\Delta_n}
\sum_{k=0}^{n-\Delta_n} \| s_k \|_{\mathsf{X}^{2},
\infty} \\
&&\quad\quad{} +
\frac{B}{N(1-\rho)^2} \sum_{k=0}^{n-1} \| s_k \|
_{\mathsf{X}^{2}, \infty}
\Biggl[ \frac{1}{\sigma_-^2} \sum_{m=1}^{k(\Delta_n)}
\frac{\| W_m(\cdot; \theta) \|_{\mathsf{X}^2, \infty}^2 \rho^{0
\vee
(k-m)}}{\{\mu g_{\theta}(Y_m)\}^2}\\
&&\hspace*{192pt}{}+ \frac{\|
W_0(\cdot; \theta) \|_{\mathsf{X}, \infty}^2 \rho^k}{\{\nu
g_{\theta}(Y_0)\}^2} \Biggr].
\end{eqnarray*}
\end{itemize}
\end{thm}

For the purpose of illustrating these bounds, assume that we are given
a set $\{ y_k ; k \geq0 \}$ of fixed observations and that all
$\| s_k \|_{\mathsf{X}^{2}, \infty}$, as well as all
fractions $\| W_k(\cdot;
\theta) \|_{\mathsf{X}^2, \infty} / \mu g_{\theta}(y_k)$, are uniformly
bounded in $k$. We then conclude that increasing the lag with $n$ as
$\Delta_n = \lceil c \log n \rceil$, $c > - 1 / \log\rho$, will
imply that $n \rho^{\Delta_n}$ tends to zero as $n$ goes to infinity,
leading to an error which is dominated by the variability due to the
particle filter (the second term of the bound in
Theorem~\ref{mainresulttheorem}(i)) and upper bounded by a quantity
proportional to
\[
\frac{1}{\sqrt{N}} \sum_{k=0}^{n-1} \Biggl[ \sum_{m=1}^{k + \lceil c
\log n
\rceil} \rho^{0 \vee(k-m)} + 1 \Biggr] \leq\frac{n}{\sqrt{N}}
\biggl(
\frac{1}{1-\rho} + 1 + \lceil
c \log n \rceil
\biggr),
\]
that is, of order $n \log n / \sqrt{N}$. Note the dependence on the mixing
coefficient $\rho$ of this rate. In contrast,
setting $\Delta_n = n$, that is, using the direct
full-path approximation, would result in a stochastic error which is
upper bounded by a quantity proportional to $n^2 / \sqrt{N}$.

%s3.2 ###
\subsection{Extension to randomly varying observations}
As mentioned, all results presented above concern smoothing
distribution approximations produced by the particle
filter algorithm \textit{conditionally} on a given sequence of observations.
In this section, we extend these results to the case of a randomly
varying observation sequence.

For the bounds presented in Theorem~\ref{mainresulttheorem}, the
conditioning on $\mathcal{G}_{n}$ can be removed by
introducing additional model assumptions. In the following, we suppose
that $\nu \ll\mu$ and that the resulting Radon--Nikod\'ym
derivative satisfies $(\mathrm{d}\nu/ \mathrm{d}\mu)_- \triangleq
\inf_{x
\in\mathsf{X}} \mathrm{d}\nu/ \mathrm{d}\mu(x) > 0$.
\begin{enumerate}
\item[(A3)]Let $t_n$ be given
by (\ref{eqtndef}). For $p \geq2$, $\ell\geq1$ and $\theta
\in\Theta$, there exists a constant $a_{p,\ell}(t_n ; \theta) \in
\mathbb{R}^+$
such that
\begin{eqnarray*}
&&\max\biggl\{ \bar{\mathbb{E}}\biggl[ \frac{\| W_k(\cdot;
\theta) \|_{\mathsf{X}^{}, \infty}^p \| s_i \|
_{\mathsf{X}^{2}, \infty}^\ell}{\{ \mu g_{\theta}(Y_k)\}^p}
\biggr], \bar{\mathbb{E}}[
\| s_i \|_{\mathsf{X}^{n+1}, \infty}^\ell] ; 0
\leq k \leq n, 0 \leq
i \leq n-1 \biggr\} \\
&&\quad\leq a_{p,\ell}(t_n ; \theta).
\end{eqnarray*}
\end{enumerate}
\begin{proposition} \label{propositionrandobs}
Assume \textup{(A1)} and
\textup{(A2)}. There then exist universal
constants $B_p$ and $B$, $B_p$ depending only on $p$, such that the
following holds true for all $N \geq1$:
\begin{itemize}
\item[(i)] if assumption \textup{(A3)} is satisfied
for $\ell= p \geq2$ and $\theta
\in\Theta$, then
\begin{eqnarray*}
\| \widehat{\gamma}^{N,\Delta_n}_{\theta, n} - \gamma
_{\theta,
n} \|_p &\leq&2 a_{p,p}^{1/p}(t_n ; \theta)
\rho^{\Delta_n} (n - \Delta_n +1) \\
&&{}+ \frac{B_p a^{1/p}_{p,p}(t_n ; \theta)}{\sqrt{N}(1-\rho)} \biggl\{
\frac{\Delta_n(n+1)}{\sigma_-} + n \biggl[ \frac{1}{\sigma
_-(1-\rho)} +
\frac{1}{( {\mathrm{d}\nu/\mathrm{d}\mu})^2_-} + 1
\biggr] \biggr\};
\end{eqnarray*}
\item[(ii)] if assumption \textup{(A3)}  is satisfied
for $p=2$, $\ell=1$ and $\theta\in\Theta$, then
\begin{eqnarray*}
| \bar{\mathbb{E}}^N_\theta [ \widehat{\gamma
}^{N,\Delta_n}_{\theta,
n} -
\gamma_{\theta, n} ] | &\leq&2 a_{2,1}(t_n ; \theta)
\rho^{\Delta_n} (n - \Delta_n +1) \\
&&{}+ \frac{B a_{2,1}(t_n ; \theta)}{N(1-\rho)^2} \Bigg\{ \frac{\Delta
_n (n+1)
}{\sigma_-^2} + n \Bigg[ \frac{1}{\sigma_-^2(1-\rho)} +
\frac{1}{({\mathrm{d}\nu/\mathrm{d}\mu} )^2_- }
\Bigg] \Bigg\}.
\end{eqnarray*}
\end{itemize}
\end{proposition}

The proof of this result is given in Section \ref{sectionrandobsproof}.
\begin{remark}
In the case of a compact state space
$\mathsf{X}$, assumption~(A3)
implies only limited
additional restrictions on the state space model. In fact, for a large
class of models, assumption~(A3)
follows as
a direct consequence of assumption~(A1).
\end{remark}

%s4 ###
\section{Applications to maximum likelihood estimation}
\label{sectionApplicationsMLE}
\setcounter{equation}{0}
We now return to the computation of the maximum likelihood estimator.
In the following, we consider models for which the set of complete data
likelihood functions is an \textit{exponential family}, that is, for
all $\theta\in\Theta$ and $n \geq0$, the joint density of
$(X_{0:n}, Y_{0:n})$ is of the form $\exp[ \langle\psi (\theta),
S_n(x_{0:n}) \rangle- c(\theta)] h(x_{0:n})$. Here, $\psi$ and the
sufficient statistics $S_n$ are
$\mathbb{R}^{d_s}$-valued functions on $\Theta$ and $\mathsf{X}^{n+1}$,
respectively, $c$ is a real-valued function on $\theta$ and $h$ is
a real-valued non-negative function on $\mathsf{X}^{n+1}$. By $\langle
\cdot, \cdot\rangle$ we denote the scalar product. All of these
functions may
depend on the observed values $y_{0:n}$, even though this
is expunged from the notation.

If the complete data likelihood function is of the particular
form above and the expectation
$\phi_{\nu, 0:n|n}(S_n; \theta)$ is finite for all $\theta
\in\Theta$, then the intermediate
quantity of EM can be written as (up to quantities which do not depend
on $\theta$) $\mathcal{Q}(\theta; \theta') = \langle\psi(\theta),
\phi_{\nu, 0:n|n}(S_n; \theta') \rangle
-c(\theta)$. Note, finally that, as mentioned in the
\hyperref[s1]{Introduction}, a typical element $S_{n,m}(x_{0:n})$, $1 \leq m
\leq d_s$, of the vector
$S_n(x_{0:n})$ is an additive functional
$S_{n,m}(x_{0:n}) = \sum_{k=0}^{n-1}s_{n,m}^{(k)}(x_{k:k+1})$
so that $\phi_{\nu, 0:n|n}(S_n; \theta')$ can be estimated
using either (\ref{eqfixedlagsmoothingestimateweighted}) or
(\ref{eqfixedlagsmoothingestimateresamp}).
Denoting by $\widehat{S}_n$ such an estimate, we may approximate the
intermediate quantity by
\[
\widehat{\mathcal{Q}}^N(\theta; \theta') = \langle\psi
(\theta),
\widehat{S}_n
\rangle-c(\theta).
\]
In the next step -- referred to as the
M-step -- $\widehat{\mathcal{Q}}^N(\theta; \theta')$ is maximized
with respect to
$\theta$, providing a new parameter estimate. This procedure is
repeated recursively given an initial guess $\widehat{\theta}_0$.

As an illustration, we consider the problem of inference in a noisily
observed AR(1) model and the stochastic volatility (SV) model. None of
these examples satisfy assumption~(A1); however,
geometric ergodicity for the models in question can be established
using bounds presented by Douc et al.
\cite{doucfortmoulinespriouret2007}. Although these bounds are
somewhat more involved than those presented in
Theorem~\ref{mixingconditionalchainsth} (e.g., the former depend on
the initial distributions and the observations), we may,
nevertheless, expect that the conclusion reached in
Section~\ref{sectiontheoretics}, that is, that the error of the
fixed-lag approximation is controlled by a lag of order $\log
n$, still applies. The situation is complicated, however,  by the fact
that the
mixing rates depend on the observations and are uniform only under the
expectation operator. In other words, there may be occasional outcomes
for which mixing is poor, even if the average
performance of the system is satisfactory.

\begin{example}[({SMCEM for noisily observed AR(1) model)}] \label{AR(1)1}
We consider the state space model
\begin{eqnarray*}
X_{k+1}&=&a X_k + \sigma_w W_{k+1},\\
Y_k&=&X_k + \sigma_v V_k
\end{eqnarray*}
with $\{W_k; k \geq1\}$ and $\{V_k; k \geq0 \}$ being
mutually independent sets of standard normal distributed variables
such that $W_{k+1}$ is independent of $(X_i, Y_i)$, $0 \leq i \leq k$,
and $V_k$ is independent of $X_k$, $(X_i, Y_i)$, $0 \leq i \leq
k-1$. The initial
distribution is chosen to be a diffuse prior so that
$\phi_{\nu,0|0}$ is $\mathcal{N}(y_0, \sigma_v^2)$. Throughout the
experiment,
we use a fixed sequence of observations produced by simulation under the
parameters $a^\ast= 0.98$,
$\sigma_w^\ast= 0.2$ and $\sigma_v^\ast= 1$. In this case,
$\psi(\theta)= [1/2 \sigma_w^2, -a/\sigma_w^2, a^2/(2
\sigma_w^2), 1/(2 \sigma_v^2)]$ and the components of the $\mathbb
{R}^4$-valued
function $x_{0:n} \mapsto S_n(x_{0:n})$ are given by
$S_{n,1}(x_{0:n})\triangleq\sum_{k=1}^{n-1} x_k^2$,
$S_{n,2}(x_{0:n}) \triangleq\sum_{k=0}^{n-1} x_k x_{k+1}$,
$S_{n,3}(x_{0:n}) \triangleq\sum_{k=0}^n x_k^2$ and
$S_{n,4}(x_{0:n}) \triangleq\sum_{k=0}^n (y_k-x_k)^2$.
Furthermore, up to terms not depending on parameters, $c(\theta) = n
\log(\sigma_w^2)/2 + (n+1) \log(\sigma_v^2)/2$. In this setting, one
step of the MCEM algorithm is carried out in the following way. Having
produced an estimate $\widehat{\theta}^{i-1}$ of the parameters
$\theta=(a, \sigma_w^2, \sigma_v^2)$ at the previous iteration, we
compute an approximation $\widehat{S}_n=(\widehat{S}_{n,1},
\widehat{S}_{n,2}, \widehat{S}_{n,3}, \widehat{S}_{n,4})$ of
$\phi_{\nu, 0:n|n}(S_n; \widehat{\theta}^{i-1})$
using the particle filter and update the parameters according to
\[
\widehat{a}^i=\frac{\widehat{S}_{n,2}}{\widehat{S}_{n,1}}, \qquad
( \widehat{\sigma}_w^i
)^2=\frac{1}{n}(\widehat{S}_{n,3} - \widehat{a}^i
\widehat{S}_{n,2}
), \qquad(\widehat{\sigma}_v^i)^2 =
\frac{\widehat{S}_{n,4}}{n+1}.
\]

We simulated, for each $n={}$100, 1000, 10,000
observations, 1000 SMC estimates of
$\phi_{\nu, 0:n|n}S_1$ using the fixed-lag smoothing technique for the
parameter values $a=0.8$, $\sigma_w=0.5$ and $\sigma_v=2$. Here, the
standard bootstrap particle filter with systematic resampling was
used, with $R_{k} \equiv Q_{}$ for all $k \geq0$. The dotted
lines indicate the exact expected values, obtained by means of
disturbance smoothing. To study the bias-variance
trade-off -- discussed in detail in the previous section -- of
the method, we used six different lags for each $n$ and
a constant particle population size $N = 1000$. The result is
displayed in Figure~\ref{figcomplags}, from which it is evident that
the bias
is controlled for a size of  lag that increases approximately
logarithmically with $n$. In particular, from the plot, we deduce that an optimal
outcome is gained when lags of size $2^4$, $2^4$ and
$2^5$ are used for $n$ being 100, 1000 and 10,000,
respectively.%\looseness=1

When the lag is sufficiently large so that we can ignore
the term of the bias which is deduced from forgetting
arguments  (being roughly of magnitude $n \rho^{\Delta_n}$),
increasing the lag further exclusively leads to an increase of
variance, as well as bias, of the estimates; compare the two last boxes
of each plot. This is completely in accordance with the
theoretical results of Section~\ref{sectionParticleApproximation}.
Note that the scale on the $y$-axis is
the same for the three panels, although the $y$-axis has been shifted in each
panel due to the fact that the value of the normalized smoothed statistic
evolves as the number of observations increases.%\looseness=1
%f2 ###
%
\begin{figure}

\includegraphics{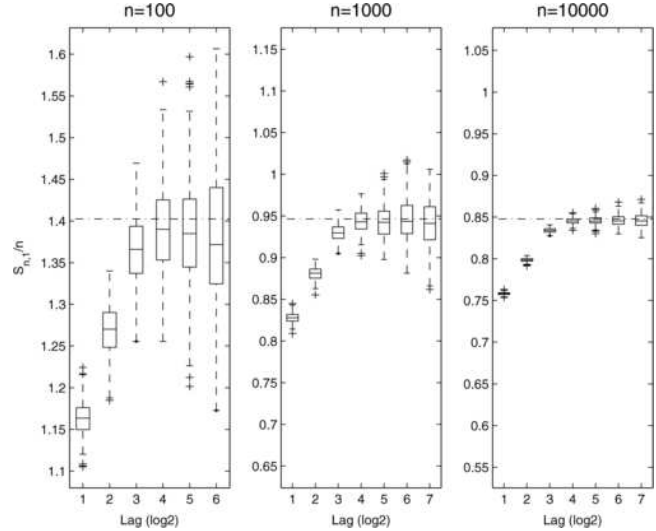}

\caption{Boxplots of estimates of $\phi_{\nu, 0:n|n}S_{n,1}/n$,
produced with the fixed-lag technique, for the noisily observed AR(1)
model in Example~\protect\ref{AR(1)1}.}
\label{figcomplags}
\end{figure}
%f3 ###
%
\begin{figure}

\includegraphics{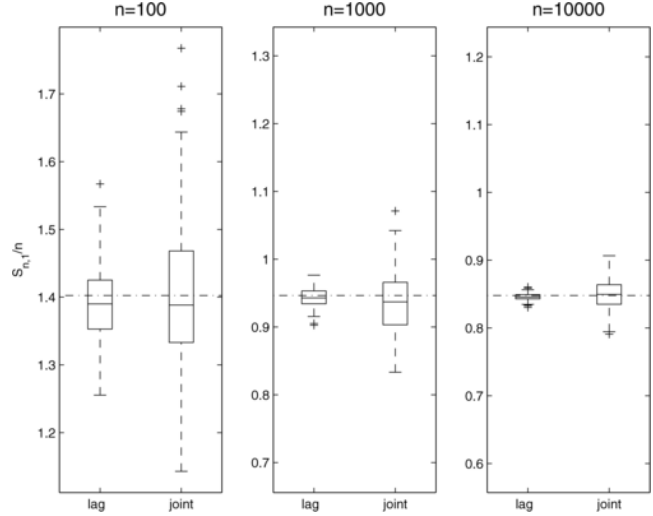}

\caption{Boxplots of estimates of $\phi_{\nu, 0:n|n}S_{n,1}/n$,
produced by means of both the fixed-lag technique and standard
trajectory-based smoothing, for the noisily observed AR(1) model in
Example~\protect\ref{AR(1)1}. Each box is based on 200 estimates,
and the size of the particle population was $N={}$1000 for all
cases.}
\label{figAR1functapprox}
\end{figure}

In Figure~\ref{figAR1functapprox}, we again report the cases
$n=100$, 1000, 10,000 observations and compare the basic
approximation strategy (\ref{equpdateapprox}) with
the one based on fixed-lag smoothing with suitable lags. Guided by the
plots of Figure~\ref{figcomplags} and the theory developed in the
previous section, we choose the lags $2^4$, $2^4$ and $2^5$,
respectively. The
number of particles was set to 1000 for all
$n$. It is obvious that fixed-lag smoothing drastically reduces the
variance without significantly raising the bias. As in the previous
figure, dotted lines indicate exact values. As expected, the
bias of the two techniques increases with $n$ since the number of
particles is held constant.\looseness=-1
\end{example}

\begin{example}[{(SMCEM for the stochastic volatility (SV) model)}]
\label{examplestovol2}
In the discrete-time case, the canonical version of the SV model
(Hull and White~\cite{hullwhite1987}, Jacquier \textit{et al.} \cite{jaquierpolsonpresnell1994}) is given by the
two relations
\begin{eqnarray*}
X_{k+1} &=& \alpha X_k + \sigma\epsilon_{k+1}, \\
Y_k &=& \beta\exp(X_k / 2) \varepsilon_k
\end{eqnarray*}
with $\{\epsilon_k; k \geq1\}$ and $\{\varepsilon_k; k \geq0 \}$ being
mutually independent sets of standard normal distributed variables
such that $W_{k+1}$ is independent of $(X_i, Y_i)$, $0 \leq i \leq k$,
and $V_k$ is independent of $X_k$, $(X_i, Y_i)$, $0 \leq i \leq
k-1$.\looseness=-1

To use the SV model in practice, we need to estimate the parameters
$\theta= (\beta, \alpha, \sigma)$. Throughout this example, we will
use a sequence of data obtained by simulation under the parameters
$\beta^\ast= 0.63$, $\alpha^\ast= 0.975$ and $\sigma^\ast=
0.16$. These
parameters are consistent with empirical estimates for daily equity
return series and are often used in simulation studies. In conformity with
Example \ref{AR(1)1}, we assume that the latent chain is initialized
by an improper diffuse prior. The SV model is within the scope of
exponential families, with $\psi(\theta) = [-\alpha^2 /(2 \sigma^2),
-1/(2 \sigma^2), \alpha/\sigma^2, -1/(2 \beta^2)]$ and components of
$S_n(x_{0:n})$ given by $S_{n,1}(x_{0:n})
\triangleq\sum_{k=0}^{n-1} x_k^2$, $S_{n,2}(x_{0:n})\triangleq
\sum_{k=1}^n x_k^2$, $S_{n,3}(x_{0:n})\triangleq\sum_{k=1}^n x_k
x_{k-1}$ and $S_{n,4}(x_{0:n}) \triangleq\sum_{k=0}^n y_k \exp(-x_k)$.
In addition, up to terms not depending on parameters, $c(\theta) =
(n+1) \log(\beta^2)/2 + (n+1) \log(\sigma^2)/2$.

Let $\widehat{S}_n=(\widehat{S}_{n,1}, \widehat{S}_{n,2},
\widehat{S}_{n,3}, \widehat{S}_{n,4})$ be a particle
approximation of $\phi_{\nu, 0:n|n}(S_n;
\widehat{\theta}^{i-1})$. To apply the Monte Carlo EM algorithm to the
SV model is not more involved than for the autoregressive model in Example
\ref{AR(1)1}. In fact, the updating formulas appear to be
completely analogous:
\[
\widehat{\alpha}^i =\frac{\widehat{S}_{n,3}}{\widehat{S}_{n,1}},
\qquad( \widehat{\sigma}^i
)^2=\frac{1}{n}(\widehat{S}_{n,2} - \widehat{\alpha}^i
\widehat{S}_{n,3}),
\qquad(\widehat{\beta}^i)^2 =\frac{\widehat
{S}_{n,4}}{n+1}.
\]

As proposal kernel $R_{k}$, we use an approximation, used by Capp\'e
\textit{et al.} (\cite{cappemoulinesryden2005}, Example~7.2.5)
and inspired
by Pitt and Shepard  \cite{pittshephard1999}, of the so-called
\textit{optimal kernel}, that is, the conditional density of $X_{k+1}$
given \textit{both} $X_k$ and $Y_{k+1}$.

%f4 ###
%
\begin{figure}

\includegraphics{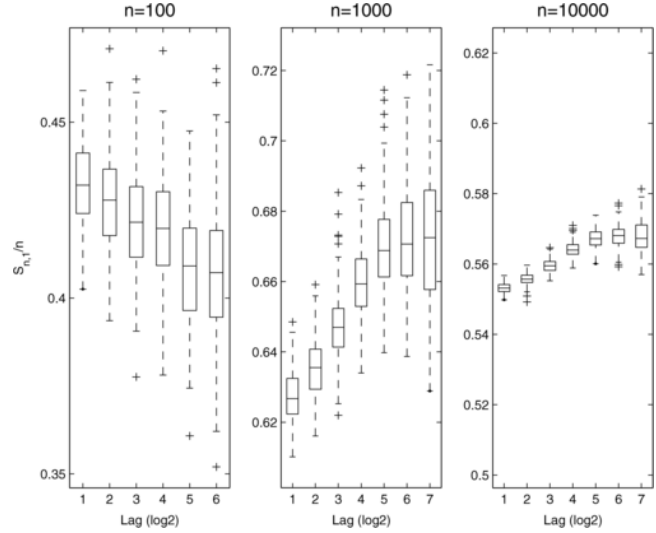}

\caption{Boxplots of estimates of $\phi_{\nu, 0:n|n}S_{n,1}/n$, produced
with the fixed-lag technique, for the SV model in
Example~\protect\ref{examplestovol2}. Each box is based on 200 estimates
and the size of the particle population was set to $N=1000$ in all
cases.}
\label{figstovolfunctapprox}
\end{figure}

We repeat the numerical investigations
of Example~\ref{AR(1)1}. The resulting approximation of
$\phi_{\nu, 0:n|n}S_{n,1}$, displayed in Figure~\ref{figstovolfunctapprox},
behaves similarly. Here, again, we observe that moderate values of the lag
$\Delta$ are sufficient to suppress the bias.

%f5 ###
%
\begin{figure}

\includegraphics{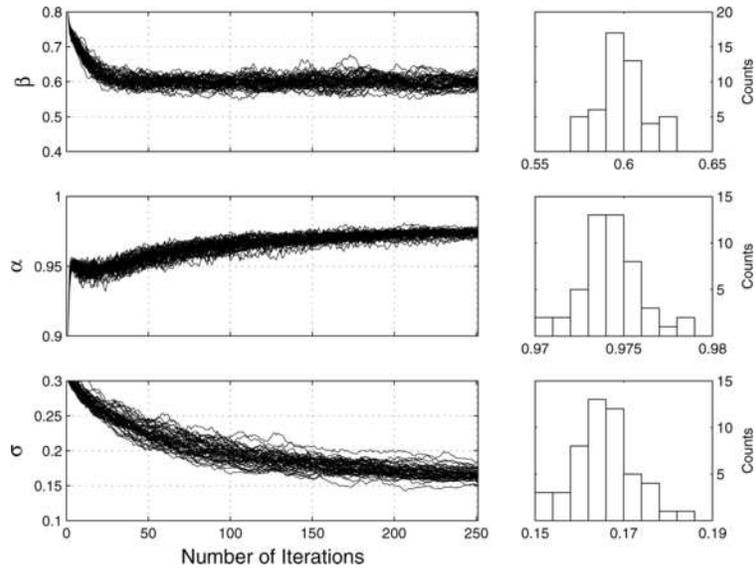}

\caption{SMCEM parameter estimates of $\beta$, $\alpha$ and $\sigma$
from $n = 5000$ observations using the standard trajectory-based
smoothing approximation. Each plot overlays 50 realizations of the
particle simulations; the histograms pertain to the final (250th)
SMCEM iteration.}
\label{figMCEMjoint}
\end{figure}

%f6 ###
%
\begin{figure}

\includegraphics{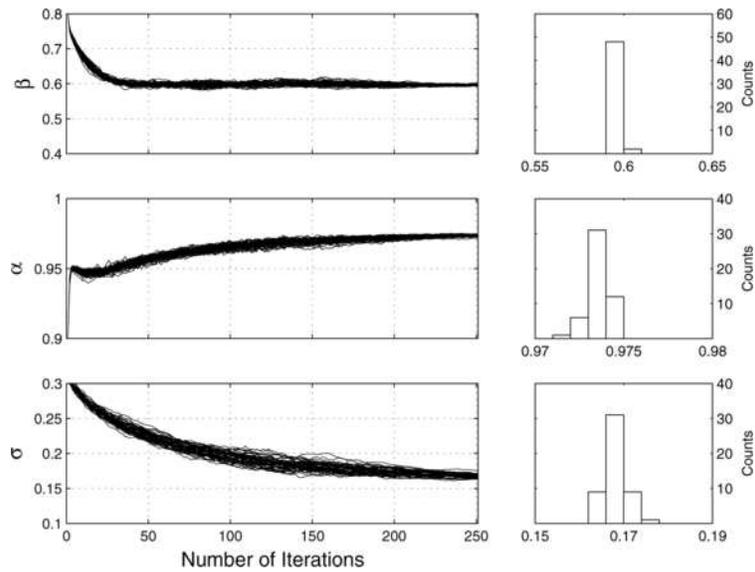}

\caption{SMCEM parameter estimates of $\beta$, $\alpha$ and $\sigma$
from $n=5000$ observations using the fixed-lag smoothing
approximation with $\Delta= 40$. Each plot overlays 50 realizations
of the particle simulations; the histograms pertain to the final
(250th) SMCEM iteration.}
\label{figMCEMlag}
\end{figure}

We finally compare the SMCEM parameter estimates obtained with the fixed-lag
approximation and with the standard trajectory-based approximation on a simulated
dataset of length $n=5000$. Note that for the SMCEM procedure to
converge to the MLE, it is necessary to increase the number of
simulations that are performed as we progress through the EM
iterations. We follow the recommendation of Fort and Moulines \cite{fortmoulines2003}
and start by running 150 iterations of the Monte Carlo EM procedure
with the number of particles set at $N=100$.
For the subsequent 100 iterations, the number of particles increases at a
quadratic rate with a final value (for the 250th Monte Carlo EM iteration)
equal to $N=1600$. The cumulative number of simulations performed
during the
250 SMCEM iterations is equal to 75,000 (times the length of the observation
sequence), which is quite moderate for a Monte Carlo-based optimization method.
In Figures~\ref{figMCEMjoint} and~\ref{figMCEMlag}, we display the
superimposed trajectories of parameter estimates for 50 realizations of the
particles, together with histograms of the final estimates (at
iteration 250)
when using, respectively, the trajectory-based approximation (in
Figure~\ref{figMCEMjoint}) and the fixed-lag approximation with
$\Delta=40$
(in Figure~\ref{figMCEMlag}).
%
%t1 ###
%
\begin{table}[b]
\caption{Mean and standard deviation of SMCEM parameter estimates at
the 250th iteration (estimated from 50 independent runs)}
\label{tabSMCEMmean_std}
\begin{tabular*}{\textwidth}{@{\extracolsep{\fill}}lccc@{}} \hline
Smoothing algorithm &
\multicolumn{1}{c}{${\hat{\beta}}$} &
\multicolumn{1}{c}{${\hat{\alpha}}$} &
\multicolumn{1}{c}{${\hat{\sigma}}$} \\\hline
Trajectory-based, & \phantom{std. }0.5991 & \phantom{std. }0.9742 & \phantom{std. }0.1659 \\
with 75,000 total simulations & std. 0.0136 & std. 0.0019 & std. 0.0070
\\
Trajectory-based, & \phantom{std. }0.5990 & \phantom{std. }0.9739 & \phantom{std. }0.1666 \\
with 750,000 total simulations & std. 0.0045 & std. 0.0011 &
std. 0.0043 \\
Fixed-lag, & \phantom{std. }0.5962 & \phantom{std. }0.9735 & \phantom{std. }0.1682 \\
with 75,000 total simulations & std. 0.0019 & std. 0.0006 & std. 0.0024
\\ \hline
\end{tabular*}
\end{table}
Not surprisingly, the fact that the particle simulations are iterated for
several successive values of the parameter estimates only amplifies the
differences observed so far. With the fixed-lag approximation, the standard
deviation of the final SMCEM parameter estimate is divided by a factor
of seven for $\beta$, and  three for $\alpha$ and $\sigma$, which
is quite
impressive in the
context of Monte Carlo methods: to achieve the same accuracy with the
trajectory-based approximation, one would need about ten times more particles
to compensate for the higher simulation variance.
Table~\ref{tabSMCEMmean_std} shows that the fixed-lag approximation (third
row) indeed remains more reliable than the trajectory-based
approximation, even
when the latter is computed from ten times more particles (second row).
Note that, for the trajectory-based approximation, multiplying the
number of
particles by ten does not reduce the standard deviation of the
estimates as
much as expected from the asymptotic theory. This is certainly due to the
moderate number of particles used in the baseline setting, as we start from
$N=100$ particles during the first SMCEM iterations and terminate with
$N=1600$.
\end{example}

\begin{appendix}
%s5 ###
\section{Proofs}
\setcounter{equation}{0}
%s5.1 ###
\subsection[Proof of Theorem 3.3]{Proof of Theorem \protect\ref{mainresulttheorem}}
The proof of Theorem \ref{mainresulttheorem} partly comprises
the geometric ergodicity of the time-reversed conditional hidden chain
(Theorem \ref{mixingconditionalchainsth}), partly
the next proposition. In the following, we omit $\theta$ from the
notation for brevity. Moreover, let $\mathcal{C}_i(\mathsf{X}^{n+1})$
be the set of bounded
measurable functions $f$ on $\mathsf{X}^{n+1}$, possibly depending on
$Y_{0:n}$, of type $f(x_{0:n}) =
\bar{f}(x_{i:n})$ for some function
$\bar{f} \dvtx \mathsf{X}^{n-i+1} \rightarrow\mathbb{R}$.
\begin{proposition} \label{propositionfixedobservations}
Assume \textup{(A1)} and \textup{(A2)}, and
let $f \in\mathcal{C}_i(\mathsf{X}^{n+1})$, $0 \leq i \leq n$.
There then exist universal
constants $B_p$ and $B$, $B_p$ depending only on $p$, such that the following
holds   for all $N \geq1$:
\begin{itemize}
\item[(i)] for all $p \geq2$,
\begin{eqnarray*}
&&\| \widehat{\phi}_{\nu, 0:n|n}^N f_i-
\phi_{\nu, 0:n|n} f_i \|_{p| \mathcal{G}_{n}} \\
&&\quad\leq\frac{B_p \| f_i \|_{\mathsf{X}^{n+1}, \infty}}{\sqrt
{N}(1-\rho)} \Biggl[
\frac{1}{\sigma_-} \sum_{k=1}^n
\frac{\| W_k \|_{\mathsf{X}^2, \infty} \rho^{0 \vee
(i-k)}}{\mu g_{}(Y_k) }+
\frac{\| W_0 \|_{\mathsf{X}, \infty} \rho^i}{\nu g_{0}} + 1
\Biggr];
\end{eqnarray*}
\item[(ii)]
\begin{eqnarray*}
&&\bigl| \bar{\mathbb{E}}^N [  \widehat{\phi}_{\nu,
0:n|n}^N f_i |
\mathcal{G}_{n}
] - \phi_{\nu, 0:n|n} f_i \bigr|\\
&&\quad\leq\frac{B \| f_i \|_{\mathsf{X}^{n+1}, \infty}}{N(1-\rho)^2}
\Biggl[
\frac{1}{\sigma_-^2} \sum_{k=1}^n
\frac{\| W_k \|_{\mathsf{X}^2, \infty}^2 \rho^{0 \vee
(i-k)}}{\{\mu g_{}(Y_k)\}^2} + \frac{\| W_0
\|_{\mathsf{X}, \infty}^2 \rho^i}{\{\nu g_{}(Y_0)\}^2} \Biggr].
\end{eqnarray*}
\end{itemize}
\label{fixed_observations}
\end{proposition}

To prove Proposition \ref{propositionfixedobservations}, we need some
preparatory lemmas and definitions. In accordance with the
mutation-selection procedure presented in Section
\ref{sectionParticleApproximation}, we have, for $k \geq1$,
$A \in\mathcal{X}^{\otimes(k+1)}$ and $i \in\{1, \ldots, N \}$, that
\begin{eqnarray*}
&&\bar{\mathbb{P}}^N
( {\xi}_{0:k}^{N,j} \in A |
\mathcal{G}_{k} \vee
\mathcal{F}_{k-1}^N) \\
&&\quad=\sum_{j=1}^N \bar{\mathbb{P}}^N
(  I^{N,i}_{k-1} = j | \mathcal{G}_{k} \vee
\mathcal{F}_{k-1}^N) \bar{\mathbb{P}}^N ( {\xi}_{0:k}^{N,j}
\in A | I^{N,i}_{k-1} = j, \mathcal{G}_{k} \vee
\mathcal{F}_{k-1}^N) \\
&&  \quad= \sum_{j=1}^N \frac{\omega_{k-1}^{N, j}}{\Omega_{k-1}^N}
R_{k-1}^\mathrm{p}(\xi_{0:k-1}^{N,j} A
).
\end{eqnarray*}
That is, conditional on $\mathcal{F}_{k-1}^N$, the
swarm $\{{\xi}_{0:k}^{N,i} ; 1 \leq i \leq N \}$ of mutated
particles at time $k$ is obtained by sampling $N$ independent and
identically distributed particles from
the measure
%
%e1 ###
%
\begin{equation} \label{eqeta}
\eta^N_k \triangleq\phi_{\nu, 0:k-1|k-1}^N R_{k-1}^\mathrm{p}.
\end{equation}
Using this, define, for $A \in\mathcal{X}^{\otimes(k+1)}$,
%
%e2 ###
%
\begin{equation} \label{eqmu}
\mu_{k|n}^N (A) \triangleq\int_A \frac{\mathrm{d}\mu
_{k|n}^N}{\mathrm{d}
\eta^N_k}(x_{0:k})  \eta^N_k(\mathrm{d}x_{0:k}),
\end{equation}
where the Radon--Nikod\'ym derivative is given by, for
$x_{0:k} \in\mathsf{X}^{k+1}$,
\[
\frac{\mathrm{d}\mu_{k|n}^N}{\mathrm{d}\eta^N_k}(x_{0:k})
\triangleq
\frac{W_k(x_{k-1:k}) L_{k} \cdots
L_{n-1}(x_{0:k}, \mathsf{X}^{n+1})}{ \phi_{\nu, 0:k-1|k-1}^N L_{k-1}
\cdots L_{n-1}(\mathsf{X}^{n+1})}.
\]
In addition, for $A \in\mathcal{X}$, let
\[
\mu_{0|n} (A) \triangleq\int_A \frac{\mu_{0|n}}{\mathrm
{d}\varsigma}(x_0)
\varsigma(\mathrm{d}x_0)
\]
with, for $x_0 \in\mathsf{X}$ and $y_0 \in\mathsf{Y}$,
\[
\frac{\mu_{0|n}}{\mathrm{d}\varsigma}(x_0) \triangleq\frac{W_0
(x_0) L_{0} \cdots L_{n-1}(x_0, \mathsf{X}^{n+1})}{
\nu[g_{}(\cdot, y_0) L_{0} \cdots L_{n-1}(\mathsf{X}^{n+1})]}.
\]
\begin{lm} \label{decompositionlemma}
Let $f \in\mathcal{B}_{\mathrm{b}}(\mathsf{X}^{n+1})$. Then, for
all $n \geq0$ and $N \geq1$,
\[
\phi_{\nu, 0:n|n}^N f - \phi_{\nu, 0:n|n} f = \sum_{k=0}^n
\varphi_k^N(f),
\]
where, for $k \geq1$,
\begin{eqnarray} \label{eqgammadef}
\varphi_k^N(f) &\triangleq&\frac{ \sum_{i=1}^N {\mathrm{d}\mu
_{k|n}^N}/{\mathrm{d}
\eta^N_k}({\xi}_{0:k}^{N,i})
\Psi_{k,n}[f]({\xi}_{0:k}^{N,i})}{\sum_{j=1}^N {\mathrm{d}
\mu_{k|n}^N}/{\mathrm{d}\eta^N_k}(\xi_{0:k}^{N,j})} - \mu_{k|n}^N
\Psi_{k,n}[f],\\
\varphi_0^N(f) &\triangleq&\frac{ \sum_{i=1}^N  {\mathrm{d}\mu
_{0|n}}/{\mathrm{d}
\varsigma}(\xi_{0}^{N,i})
\Psi_{0,n}[f](\xi_{0}^{N,i})}{\sum_{j=1}^N {\mathrm{d}
\mu_{0|n}}/{\mathrm{d}\varsigma}(\xi_{0}^{N,j})} - \mu_{0|n}
\Psi_{0,n}[f] \nonumber
\end{eqnarray}
and the operators $\Psi_{k,n} \dvtx \mathcal{B}_{\mathrm{b}}(\mathsf
{X}^{k+1}) \rightarrow\mathcal{B}_{\mathrm{b}}(\mathsf{X}^{k+1})$,
$0 \leq k \leq n+1$, are, for fixed points $\widehat{x}_{0:k} \in
\mathsf{X}^{k+1}$, defined by
%e3 ###
%
\begin{equation} \label{eqfhat}
\Psi_{k,n}[f](x_{0:k}) \triangleq\frac{L_{k} \cdots L_{n-1}
f(x_{0:k})}{L_{k} \cdots L_{n-1}(x_{0:k}, \mathsf{X}^{n+1})}-
\frac{L_{k} \cdots L_{n-1} f(\widehat{x}_{0:k})}{L_{k} \cdots
L_{n-1}(\widehat{x}_{0:k}, \mathsf{X}^{n+1})}.
\end{equation}
\end{lm}
\begin{pf}
As a starting point, consider the decomposition
\begin{eqnarray*}\label{decomposition}
&&\phi_{\nu, 0:n|n}^N f - \phi_{\nu, 0:n|n} f\\
 &&\quad=
\frac{\phi_{\nu,0}^N L_{0} \cdots L_{n-1}
f}{\phi_{\nu,0}^N
L_{0} \cdots L_{n-1}(\mathsf{X}^{n+1})} - \phi_{\nu, 0:n|n} f \\
&&\qquad{}+ \sum_{k=1}^n \biggl[ \frac{\phi_{\nu, 0:k|k}^N L_{k} \cdots
L_{n-1} f}{\phi_{\nu, 0:k|k}^N L_{k} \cdots L_{n-1}(\mathsf{X}^{n+1})}
- \frac{\phi_{\nu, 0:k-1|k-1}^N L_{k-1} \cdots L_{n-1}
f}{\phi_{\nu, 0:k-1|k-1}^N L_{k-1} \cdots L_{n-1}(\mathsf{X}^{n+1})}
\biggr].
\end{eqnarray*}
Using the definitions (\ref{eqeta}) and (\ref{eqmu}) of $\eta^N_k$
and $\mu_{k|n}^N$, respectively, we may write, for $k \geq1$,
\begin{eqnarray*}
&& \frac{\phi_{\nu, 0:k-1|k-1}^N L_{k-1} \cdots L_{n-1}
f}{\phi_{\nu, 0:k-1|k-1}^N L_{k-1} \cdots L_{n-1}(\mathsf
{X}^{n+1})} \\
&&\quad = \eta^N_k \biggl[ \frac{
W_k(\cdot) L_{k} \cdots
L_{n-1}f(\cdot)}{\phi_{\nu, 0:k-1|k-1}^N L_{k-1} \cdots
L_{n-1}(\mathsf{X}^{n+1})} \biggr] \\
&&\quad = \eta^N_k \biggl[ \frac{
W_k(\cdot) L_{k} \cdots
L_{n-1} (\cdot, \mathsf{X}^{n+1})}{\phi_{\nu, 0:k-1|k-1}^N L_{k-1}
\cdots
L_{n-1}(\mathsf{X}^{n+1})} \biggl\{ \Psi_{k,n}[f](\cdot) +
\frac{L_{k} \cdots L_{n-1} f(\widehat{x}_{0:k})}{L_{k}
\cdots L_{n-1}(\widehat{x}_{0:k}, \mathsf{X}^{n+1})} \biggr\}
\biggr] \\
&&\quad = \mu_{k|n}^N \biggl[ \Psi_{k,n}[f](\cdot) + \frac{L_{k}
\cdots L_{n-1} f(\widehat{x}_{0:k})}{L_{k} \cdots
L_{n-1}(\widehat{x}_{0:k}, \mathsf{X}^{n+1})} \biggr] \\
&&\quad = \mu_{k|n}^N \Psi_{k,n}[f] + \frac{L_{k}
\cdots L_{n-1} f(\widehat{x}_{0:k})}{L_{k} \cdots
L_{n-1}(\widehat{x}_{0:k}, \mathsf{X}^{n+1})}.
\end{eqnarray*}
On the other hand,
\begin{eqnarray*}
&&\frac{\phi_{\nu, 0:k|k}^N L_{k} \cdots L_{n-1}
f}{\phi_{\nu, 0:k|k}^N L_{k} \cdots L_{n-1}(\mathsf{X}^{n+1})}\\
&&\quad=
\frac{ \sum_{i=1}^N{\mathrm{d}\mu_{k|n}^N}/{\mathrm{d}
\eta^N_k}({\xi}_{0:k}^{N,i})
\Psi_{k,n}[f]({\xi}_{0:k}^{N,i})}{\sum_{j=1}^N {\mathrm{d}
\mu_{k|n}^N}/{\mathrm{d}\eta^N_k}(\xi_{0:k}^{N,j})} + \frac{L_{k}
\cdots L_{n-1} f(\widehat{x}_{0:k})}{L_{k} \cdots
L_{n-1}(\widehat{x}_{0:k}, \mathsf{X}^{n+1})}
\end{eqnarray*}
and, by combining the two latter identities, it follows from the
definition (\ref{eqgammadef}) of $\varphi_k^N(f)$ that, for $k
\geq1$,
\[
\varphi_k^N(f) = \frac{\phi_{\nu, 0:k|k}^N L_{k} \cdots
L_{n-1} f}{\phi_{\nu, 0:k|k}^N L_{k} \cdots L_{n-1}(\mathsf{X}^{n+1})}
- \frac{\phi_{\nu, 0:k-1|k-1}^N L_{k-1} \cdots L_{n-1}
f}{\phi_{\nu, 0:k-1|k-1}^N L_{k-1} \cdots L_{n-1}(\mathsf{X}^{n+1})}
.
\]

The identity
\[
\varphi_0^N(f) = \frac{\phi_{\nu,0}^N L_{0} \cdots L_{n-1}
f}{\phi_{\nu,0}^N
L_{0} \cdots L_{n-1}(\mathsf{X}^{n+1})} - \phi_{\nu, 0:n|n} f
\]
can be verified in a similar manner.
\end{pf}

Note that, conditional on $\mathcal{F}_{k-1}^N$, the first term on the
right-hand side
of~(\ref{eqgammadef}) is nothing but an importance sampling estimate
of $\mu_{k|n}^N \Psi_{k,n}[f]$, based on $N$ independent
$\eta^N_k$-distributed variables.
\begin{lm} \label{lemmafhatbound}
Assume \textup{(A1)} and let, for $n \geq0$ and $0
\leq i \leq
n$, $f_i \in\mathcal{C}_i(\mathsf{X}^{n+1})$. Furthermore, let, for
$k \geq
0$, the operator $\Psi_{k, n}$ be defined
via \textup{(\ref{eqfhat})}. Then,
\[
\| \Psi_{k,n}[f_i] \|_{\mathsf{X}^{k+1}, \infty} \leq2
\rho^{0 \vee(i-k)}
\| f_i \|_{\mathsf{X}^{n+1}, \infty}.
\]
\end{lm}

\begin{pf}
For $k \geq i$, we bound $\Psi_{k,n}[f_i]$ from above by $2
\| f_i \|_{\mathsf{X}^{n+1}, \infty}$; however, for $k <
i$, a geometrically decreasing
bound of the function can be obtained by using the forgetting property
of the conditional latent
chain. Hence, by the Markov property of the posterior chain and using
the definition of the forward kernels (see (\ref{eqforwardkernel})),
\begin{eqnarray*}
\frac{L_{k} \cdots L_{n-1} f_i(x_{0:k})}{L_{k} \cdots
L_{n-1}(x_{0:k}, \mathsf{X}^{n+1})}
&=& \mathbb{E}[  f_i(X_{i:n}) | X_{0:k} =
x_{0:k}, \mathcal{G}_{n}
]\\
&=& \mathbb{E}\bigl[
\mathbb{E}[  f_i(X_{i:n}) | X_i = x_i,
\mathcal{G}_{n} ] | X_k = x_k, \mathcal{G}_{n}
\bigr]\\
&=& \mathrm{F}_{k|n} \cdots\mathrm{F}_{i-1|n} \{ x_k, \mathbb
{E}[
 f_i(X_{i:n}) | X_i = \cdot, \mathcal{G}_{n} ]
\}
\end{eqnarray*}
with $x_{0:k} \in\mathsf{X}^{k+1}$. Therefore, we may, for $k < i$,
rewrite $\Psi_{k,n}[f_i](x_{0:k})$ as
\begin{eqnarray*}
&&\Psi_{k,n}[f_i](x_{0:k}) \\
&&\quad= \int_\mathsf{X}\{ \mathrm
{F}_{k|n} \cdots
\mathrm{F}_{i-1|n}(x_k, \mathrm{d}x_i) - \mathrm{F}_{k|n} \cdots
\mathrm{F}_{i-1|n}(\hat{x}_k, \mathrm{d}
x_i) \} \mathbb{E}[
 f_i(X_{i:n}) | X_i = x_i, \mathcal{G}_{n} ].
\end{eqnarray*}
Applying Theorem~\ref{mixingconditionalchainsth} to this difference
yields
\begin{eqnarray*}
&&| \Psi_{k,n}[f_i](x_{0:k}) |\\
&&\quad\leq2 \| \mathbb{E}[  f_i(X_{i:n}) | X_i =
\cdot, \mathcal{G}_{n} ] \|_{\mathsf{X}^{}, \infty}
\| \mathrm{F}_{k|n} \cdots\mathrm{F}_{i-1|n} (x_{k}, \cdot
) - \mathrm{F}_{k|n} \cdots
\mathrm{F}_{i-1|n} (\widehat{x}_{k}, \cdot) \|_{\mathrm
{TV}}\\
&&\quad\leq2 \| \mathbb{E}[  f_i(X_{i:n}) | X_i =
\cdot, \mathcal{G}_{n} ] \|_{\mathsf{X}^{}, \infty}
\rho^{i-k} \leq2 \| f_i \|_{\mathsf{X}^{n+1}, \infty}
\rho^{i-k}
.
\end{eqnarray*}
\end{pf}
\begin{lm} \label{lemmaradderbounds}
Assume \textup{(A1)} and let $n \geq0$. Then, for
all $1 \leq
k \leq
n$, $x_{0:k} \in\mathsf{X}^{k+1}$, $y_k \in\mathsf{Y}$ and $N \geq1$,
\[
\frac{\mathrm{d}\mu_{k|n}^N}{\mathrm{d}\eta^N_k}(x_{0:k}) \leq
\frac{\| W_k \|_{\mathsf{X}^{2}, \infty}}{\mu g_{}(y_k)
(1-\rho) \sigma_-},
\]
where $\eta^N_k$ and $\mu_{k|n}^N$ are defined in \textup{(\ref{eqeta})} and
\textup{(\ref{eqmu})}, respectively.
\end{lm}

\begin{pf}
First write, for $x_{0:k} \in\mathsf{X}^{k+1}$ and $y_{k+1} \in
\mathsf{Y}$,
\begin{eqnarray}\label{eqderivativenumerator}
&& L_{k} \cdots L_{n-1}(x_{0:k}, \mathsf{X}^{n+1})
\nonumber\\
&&\quad = \int_\mathsf{X}q(x_k, x_{k+1}) L_{k+1} \cdots L_{n-1}(x_{0:k+1},
\mathsf{X}^{n+1}) g_{}(x_{k+1}, y_{k+1})  \mu(\mathrm
{d}x_{k+1})\nonumber\\
&&\quad \leq\sigma_+\int_\mathsf{X}L_{k+1} \cdots L_{n-1}(x_{0:k+1},
\mathsf{X}^{n+1}) g_{}(x_{k+1}, y_{k+1})  \mu(\mathrm{d}
x_{k+1}).
\end{eqnarray}
Now, since the function $L_{k+1} \cdots L_{n-1}(\cdot,\mathsf{X}^{n+1})$
is constant in all but the last component of the argument,
\begin{eqnarray}\label{eqderivativedenominator}
&& L_{k-1} \cdots L_{n-1}(x_{0:k-1}, \mathsf{X}^{n+1})
\nonumber\\
&&\quad = \int_\mathsf{X}q(x_{k-1},x_k) g_{}(x_k, y_k) \int_\mathsf
{X}q(x_k, x_{k+1})
L_{k+1} \cdots L_{n-1}(x_{0:k+1},\mathsf{X}^{n+1}) \nonumber\\
&&\quad\hspace*{112pt}{} \times g_{}(x_{k+1}, y_{k+1})  \mu^{\otimes
2}(\mathrm{d}
x_{k:k+1}) \nonumber\\
&&\quad \geq\mu g_{}(y_k) \sigma_-^2 \int_\mathsf{X}L_{k+1} \cdots
L_{n-1}(x_{0:k+1}, \mathsf{X}^{n+1}) g_{}(x_{k+1}, y_{k+1})\mu
(\mathrm{d}
x_{k+1}).
\end{eqnarray}
Since the integrals in (\ref{eqderivativenumerator}) and
(\ref{eqderivativedenominator}) are equal, the bound of the lemma
follows.
\end{pf}

\begin{pf*}{Proof of Proposition \protect\ref{propositionfixedobservations}}
We start with (i). Since, conditional on $\mathcal{F}_{n}^N$, the random
variables $f_i(\hat{\xi}_{0:n}^{N,j})$, $1 \leq j \leq N$,
are independent and identically distributed with expectations
%
%e4 ###
%
\begin{equation} \label{eqselectionexpectation}
\bar{\mathbb{E}}^N_\theta
 [  f_i(\hat{\xi}_{0:n}^{N,j}) |
\mathcal{G}_{n}
\vee
\mathcal{F}_{n}^N]
= \frac{1}{\Omega_{n}^N} \sum_{j=1}^N
\omega_{n}^{N, j} f_i(\xi_{0:n}^{N,j}),
\end{equation}
applying the Marcinkiewicz--Zygmund inequality provides the bound
%
%e5 ###
%
\begin{equation} \label{eqweightningbound}
N^{p/2} \bar{\mathbb{E}}^N_\theta
 \Biggl[  \Bigg| \frac{1}{N} \sum_{j=1}^N
f_i(\hat{\xi}_{0:n}^{N,j})-\frac{1}{\Omega_{n}^N}
\sum_{j=1}^{N} \omega_{n}^{N, j} f_i(\xi_{0:n}^{N,j})\Bigg|^p
\Big| \mathcal{G}_{n} \vee\mathcal{F}_{n}^N\Biggr] \leq C_p
\| f_i \|_{\mathsf{X}^{n+1}, \infty}^p,
\end{equation}
where $C_p$ is a universal constant depending only  on $p$. Having
control of this discrepancy, we focus instead on the $\mathsf{L}^{p}$ error
associated with the weighted empirical measure
$\phi_{\nu, 0:n|n}^N$. We
make use of the identity
\[
a/b-c = (a/b)(1-b)+a-c
\]
on each term of the decomposition provided by Lemma
\ref{decompositionlemma}. This, together with Minkowski's inequality,
gives us the bound
%
%e7 ###
%e6 ###
%
\begin{eqnarray}\label{eqtermwisedecomposition}
\| \varphi_k^N (f_i) \|_{p |
\mathcal{G}_{n} \vee\mathcal{F}_{k}^N} &\leq&\Biggl\| \frac{1}{N}
\sum_{j=1}^N \frac{\mathrm{d}\mu_{k|n}^N}{\mathrm{d}
\eta^N_k}(\xi
_{0:k}^{N,j}) \Psi_{k,n}[f_i](\xi_{0:k}^{N,j})(k) -
\mu_{k|n}^N \Psi_{k,n}[f_i] \Biggr\|_{p |
\mathcal{G}_{n} \vee\mathcal{F}_{k-1}^N}
\nonumber\\
&&{}+ \| \Psi_{k,n}[f_i] \|_{\mathsf{X}^{k+1}, \infty}
\Biggl\|
\frac{1}{N} \sum_{j=1}^N \frac{\mathrm{d}\mu_{k|n}^N}{\mathrm{d}
\eta^{N}_k}(\xi
_{0:k}^{N,j}) -1 \Biggr\|_{p |
\mathcal{G}_{n} \vee\mathcal{F}_{k-1}^N}.
\end{eqnarray}
Applying the Marcinkiewicz--Zygmund inequality to the first term of
this bound gives
%
%e9 ###
%e8 ###
%
\begin{eqnarray} \label{eqfirstpartbound}
&&N^{p/2} \bar{\mathbb{E}}^N
\Biggl[  \Bigg| \frac{1}{N} \sum
_{j=1}^N \frac
{\mathrm{d}
\mu_{k|n}^N}{\mathrm{d}\eta^N_k}(\xi_{0:k}^{N,j})
\Psi_{k,n}[f_i](\xi
_{0:k}^{N,j})- \mu_{k|n}^N \Psi_{k,n}[f_i]
\Bigg|^p \Big| \mathcal{G}_{n} \vee\mathcal{F}_{k-1}^N \Biggr]\nonumber\\
&&\quad\leq C_p \biggl\| \frac{\mathrm{d}\mu_{k|n}^N}{\mathrm{d}\eta^N_k}
\biggr\|_{\mathsf{X}^{k+1}, \infty}^p
\| \Psi_{k,n}[f_i] \|_{\mathsf{X}^{k+1}, \infty}^p
\end{eqnarray}
and treating the second term in a similar manner yields
%
%e10 ###
%
\begin{equation} \label{eqsecondtermbound}
N^{p/2} \bar{\mathbb{E}}^N
 \Biggl[  \Bigg| \frac{1}{N} \sum
_{j=1}^N \frac
{\mathrm{d}
\mu_{k|n}^N}{\mathrm{d}\eta^N_k}(\xi_{0:k}^{N,j}) -1 \Bigg|^p
\Big| \mathcal{G}_{n} \vee\mathcal{F}_{k-1}^N \Biggr] \leq C_p
\biggl\| \frac{\mathrm{d}\mu_{k|n}^N}{\mathrm{d}\eta^N_k} \biggr\|
_{\mathsf{X}^{k+1}, \infty}^p.
\end{equation}
Thus, we obtain, by inserting these bounds into
(\ref{eqtermwisedecomposition}) and applying
Lemmas~\ref{lemmafhatbound} and \ref{lemmaradderbounds},
%
%e11 ###
%
\begin{equation} \label{eqgammabound}
\sqrt{N} \| \varphi_k^N(f_i) \|_{p | \mathcal{G}_{n}
\vee\mathcal{F}_{k-1}^N} \leq4 C_p^{1/p} \rho^{0 \vee(i-k)} \frac{
\| W_k \|_{\mathsf{X}^{2}, \infty} \| f_i \|
_{\mathsf{X}^{n+1}, \infty}}{
\mu g_{}(y_k) (1-\rho) \sigma_-}.
\end{equation}

For the first term of the decomposition provided by Lemma
(\ref{decompositionlemma}), we have, using the same decomposition
technique as in (\ref{eqtermwisedecomposition}) and repeating the
arguments of Lemma~\ref{lemmaradderbounds},
%
%e12 ###
%

%
\begin{eqnarray}\label{eqfirsttermbound}
\sqrt{N} \| \varphi_0^N(f_i) \|_{p | \mathcal{G}_{n}}
&\leq&2
C_p^{1/p} \biggl\| \frac{\mu_{0|n}}{\mathrm{d}\varsigma} \biggr\|
_{\mathsf{X}^{}, \infty}
\| \Psi_{0:n}[f_i] \|_{\mathsf{X}^{}, \infty} \nonumber\\[-8pt]\\[-8pt]
\nonumber&\leq&4 C_p^{1/p} \rho^i \frac{ \| W_0 \|_{\mathsf
{X}^{}, \infty}
\| f_i \|_{\mathsf{X}^{n+1}, \infty}}{ \nu g_{}(y_0)
(1-\rho)}.
\end{eqnarray}
Now, (i) follows by a straightforward application of Minkowski's
inequality together with (\ref{eqweightningbound}),
(\ref{eqgammabound}) and (\ref{eqfirsttermbound}).

We turn to (ii). By means of the identity
\[
a/b-c = (a/b)(1-b)^2+(a-c)(1-b)+c(1-b)+a-c
\]
applied to~(\ref{eqgammadef}), we obtain the bound
\begin{eqnarray*}
&&\bigl| \bar{\mathbb{E}}^N [  \varphi_k^N(f_i) |
\mathcal{G}_{n}
\vee\mathcal{F}_{k-1}^N ] \bigr|\\
&&\quad\leq\| \Psi_{k,n}[f_i] \|_{\mathsf{X}^{k+1}, \infty}
\Biggl\| \frac{1}{N}
\sum_{j=1}^N \frac{\mathrm{d}\mu_{k|n}^N}{\mathrm{d}\eta
^{N}_{k}}(\xi_{0:k}^{N,j})
-1 \Biggr\|_{2 | \mathcal{G}_{n} \vee\mathcal{F}_{k-1}^N}^2 \\
&&\quad\quad{}+ \Biggl\| \frac{1}{N} \sum_{j=1}^N \frac{\mathrm{d}\mu
_{k|n}^N}{\mathrm{d}
\eta^N_k}(\xi
_{0:k}^{N,j}) \Psi_{k,n}[f_i](\xi_{0:k}^{N,j}) -
\mu_{k|n}^N \Psi_{k,n}[f_i] \Biggr\|_{2| \mathcal{G}_{n} \vee
\mathcal{F}_{k-1}^N} \\
&&\quad\quad\quad{}\times\Biggl\| \frac{1}{N} \sum_{j=1}^N \frac{\mathrm{d}
\mu_{k|n}^N}{\mathrm{d}\eta^N_k}(\xi_{0:k}^{N,j}) - 1
\Biggr\|_{2| \mathcal{G}_{n} \vee\mathcal{F}_{k-1}^N}.
\end{eqnarray*}
Thus, we get, by reusing (\ref{eqfirstpartbound}) and
(\ref{eqsecondtermbound}),
\begin{eqnarray}\label{eqgammabiasbound}
\bigl| {\bar{\mathbb{E}}^N} [  \varphi^N_k(f_i) |
\mathcal{G}_{n}
] \bigr|
&\leq&{\bar{\mathbb{E}}^N} \bigl[ \bigl| \bar{\mathbb{E}}^N_\theta
 [
 \varphi_k^N(f_i) | \mathcal{G}_{n} \vee
\mathcal{F}_{k-1}^N ] \bigr| | \mathcal{G}_{n} \bigr]
\nonumber\\[-8pt]\\[-8pt]\nonumber
&\leq&4 C_2 \rho^{0 \vee(i-k)} \frac{\| W_k \|_{\mathsf
{X}^{2}, \infty}^2
\| f_i \|_{\mathsf{X}^{n+1}, \infty}}{N \{\mu g_{}(y_k)\}
^2 (1 - \rho)^2 \sigma_-^2}
\end{eqnarray}
and treating the last term of the decomposition in a completely
similar manner yields
%
%e13 ###
%
\begin{equation} \label{eqfirsttermbiasbound}
\bigl| {\bar{\mathbb{E}}^N} [  \varphi^N_0(f_i) |
\mathcal{G}_{n}
] \bigr| \leq4 C_2 \rho^i \frac{\| W_0 \|
_{\mathsf{X}^{}, \infty}^2
\| f_i \|_{\mathsf{X}^{n+1}, \infty}}{N \{ \nu
g_{}(y_0) \}^2 (1 - \rho)^2}.
\end{equation}
Finally, from (\ref{eqselectionexpectation}), we conclude that the
multinomial selection mechanism does not introduce any additional
bias and, consequently, (ii) follows from the
triangle inequality, together with (\ref{eqgammabiasbound}) and
(\ref{eqfirsttermbiasbound}).
\end{pf*}

Having established Proposition
\ref{propositionfixedobservations}, we are now ready to proceed
with the proof of Theorem~\ref{mainresulttheorem}.

\begin{pf*}{Proof of Theorem \protect\ref{mainresulttheorem}}
Decomposing the difference in question yields the bound
%
%e15 ###
%e14 ###
%
\begin{eqnarray} \label{eqsumdecomposition}
\| \hat{\gamma}^{N,\Delta_n}_n - \gamma_n \|_{p |
\mathcal{G}_{n}}
&\leq&\sum_{k=0}^{n-1} \bigl\| \widehat{\phi}_{\nu, 0:k(\Delta
_n)|k(\Delta_n)}^N s_k -
\phi_{\nu, 0:k(\Delta_n)|k(\Delta_n)} s_k \bigr\|_{p|\mathcal
{G}_{n}} \nonumber\\[-8pt]\\[-8pt]\nonumber
&&{}+ \sum_{k=0}^{n - \Delta_n} | \phi_{\nu, 0:k + \Delta_n|k +
\Delta_n} s_k
-\phi_{\nu, 0:n|n}s_k |,
\end{eqnarray}
where we have set $k(\Delta_n)= (k+\Delta_n) \wedge n$.
By writing
\begin{eqnarray*}
&& \mathbb{E}[  s_k(X_k, X_{k+1}) | X_{k+\Delta_n+1}
= x_{k +\Delta_n+1},
\mathcal{G}_{k + \Delta_n} ] \\
&&\quad = \mathbb{E}\bigl[  \mathbb{E}[  s_k(X_k,
X_{k+1}) | X_{k+1}
= x_{k+1},
\mathcal{G}_{k+\Delta_n} ] | X_{k+\Delta_n+1} =
x_{k+\Delta_n+1},
\mathcal{G}_{k+\Delta_n} \bigr] \\
&&\quad = \mathrm{B}_{\nu, k+\Delta_n|k+\Delta_n} \cdots\mathrm{B}_{\nu,
k+1|k+\Delta_n} ( x_{k+\Delta_n+1},
\widehat{s}_{k|k+\Delta_n})
\end{eqnarray*}
with, for $x \in\mathsf{X}$,
\[
\widehat{s}_{k|k+\Delta_n}(x) \triangleq\mathbb{E}[
s_k(X_k, X_{k+1})
| X_{k+1} = x, \mathcal{G}_{k+\Delta_n} ],
\]
we get that
\begin{eqnarray*}
&&\phi_{\nu, 0:k + \Delta_n|k + \Delta_n} s_k
- \phi_{\nu, 0:n|n} s_k \\
&&\quad = \psi_{k+\Delta_n+1|k+\Delta_n} \mathrm{B}_{\nu,
k+\Delta_n|k+\Delta_n} \cdots
\mathrm{B}_{\nu,
k+1|k+\Delta_n} (\widehat{s}_{k|k+\Delta_n}) \\
&&\quad\quad{}  - \psi_{k+\Delta_n+1|n} \mathrm{B}_{\nu,
k+\Delta_n|k+\Delta_n} \cdots
\mathrm{B}_{\nu, k+1|k+\Delta_n} ( \widehat{s}_{k|k+\Delta_n}
),
\end{eqnarray*}
where we have defined, for $\ell, m \geq0$, $\psi_{\ell|m}
\triangleq
\operatorname{\mathbb{P}}( X_\ell\in\cdot| \mathcal{G}_{m})$.
Hence, we obtain, using the
exponential forgetting property (see
Theorem~\ref{mixingconditionalchainsth}) of the time-reversed
conditional hidden chain,
\begin{eqnarray}\label{eqbackwardmixing}
&& | \phi_{\nu, 0:k + \Delta_n|k + \Delta_n} s_k
- \phi_{\nu, 0:n|n} s_k | \nonumber\\
&&\quad \leq2 \| \widehat{s}_{k|k+\Delta_n} \|_{\mathsf
{X}^{}, \infty} \|
\psi_{k+\Delta_n+1|k+\Delta_n} \mathrm{B}_{\nu,
k+\Delta_n|k+\Delta_n} \cdots\mathrm{B}_{\nu,
k+1|k+\Delta_n} ( \cdot)  \nonumber\\[-8pt]\\[-8pt]
&&\quad\quad{}  - \psi_{k+\Delta_n+1|n} \mathrm{B}_{\nu,
k+\Delta_n|k+\Delta_n} \cdots\mathrm{B}_{\nu,
k+1|k+\Delta_n} ( \cdot)
\|_{\mathrm{TV}} \nonumber\\\nonumber
&&\quad \leq2 \rho^{\Delta_n} \| s_k \|_{\mathsf{X}^{2},
\infty}.
\end{eqnarray}
Substituting (\ref{eqbackwardmixing}) and the bound of Proposition
\ref{propositionfixedobservations}(i) into the
decomposition (\ref{eqsumdecomposition}) completes the proof of (i).
The proof of part (ii) is entirely analogous and is omitted
for brevity.
\end{pf*}
%
%s5.2 ###
\subsection[Proof of Proposition 3.4]{Proof of Proposition \protect\ref{propositionrandobs}}
\label{sectionrandobsproof}
\begin{enumerate}
\item[(A4)]Let $f_i$ be the function of Proposition
\ref{propositionfixedobservations}. For $p \geq2$, $\ell\geq1$,
there exists a constant $\alpha_{p,\ell}^{(n)}(f_i) \in\mathbb
{R}^+$ such that
\[
\max\biggl\{ \bar{\mathbb{E}}\biggl[ \frac{\| W_k \|
_{\mathsf{X}^{}, \infty}^p
\| f_i \|_{\mathsf{X}^{n+1}, \infty}^\ell}{\{\mu
g_{}(Y_k)\}^p} \biggr], \bar{\mathbb{E}}[ \| f_i
\|_{\mathsf{X}^{n+1}, \infty}^\ell
] ; 0 \leq k
\leq n \biggr\} \leq\alpha_{p,\ell}^{(n)}(f_i).
\]
\end{enumerate}

Under assumption~(A4), we have the
following result.
\begin{proposition} \label{propositionrandobsfi}
Assume \textup{(A1)} and
\textup{(A2)}. There then exist universal
constants $B_p$ and $B$, $B_p$ depending only on $p$, such that the
following holds true for all $N \geq1$:
\begin{itemize}
\item[(i)] if assumption \textup{ (A4)}
is satisfied
for $\ell= p \geq2$, then
\[
\| \widehat{\phi}_{\nu, 0:n|n}^N f_i-
\phi_{\nu, 0:n|n} f_i \|_p \leq
\frac{B_p [\alpha_{p,p}^{(n)}(f_i)]^{1/p}}{\sqrt{N}(1-\rho)} \biggl[
\frac{1-\rho^i}{\sigma_-(1-\rho)} + \frac{n-i}{\sigma_-} +
\frac{\rho^i}{( {\mathrm{d}\nu}/{\mathrm{d}\mu}
)_-} + 1
\biggr] ;
\]
\item[(ii)] if assumption \textup{(A4)} is satisfied
for $p=2$, $\ell= 1$, then
\[
| {\bar{\mathbb{E}}^N} [ \widehat{\phi}_{\nu, 0:n|n}^N f_i -
\phi_{\nu,
0:n|n} f_i ] | \leq
\frac{B \alpha_{2,1}^{(n)}(f_i)}{N(1-\rho)^2} \bigg[
\frac{1-\rho^i}{\sigma_-^2(1-\rho)} + \frac{n-i}{\sigma_-^2} +
\frac{\rho^i}{( {\mathrm{d}\nu}/{\mathrm{d}\mu}
)^2_-} \bigg].
\]
\end{itemize}
\end{proposition}

\begin{pf}
The proof of the first part is straightforward: combining Proposition
\ref{propositionfixedobservations} and Minkowski's inequality
provides the bound
\begin{eqnarray*}
&& \| \widehat{\phi}_{\nu, 0:n|n}^N f_i-
\phi_{\nu, 0:n|n} f_i \|_p\\
&&\quad = \bar{\mathbb{E}}^{1/p} [ \| \widehat{\phi}_{\nu,
0:n|n}^N f_i-
\phi_{\nu, 0:n|n} f_i \|_{p | \mathcal{G}_{n}}^p ]\\
&&\quad \leq\frac{B_p }{\sqrt{N}(1-\rho)} \Biggl\{
\frac{1}{\sigma_-} \sum_{k=1}^n \bar{\mathbb{E}}^{1/p} \biggl[
\frac{\| W_k \|_{\mathsf{X}^{}, \infty}^p \| f_i
\|_{\mathsf{X}^{n+1}, \infty}^p}{\{\mu
g_{}(Y_k)\}^p} \biggr] \rho^{0 \vee(i-k)}  \\
&&\hspace*{64pt}\quad{} +  \frac{1}{( {\mathrm{d}
\nu}/{\mathrm{d}\mu})_-}
\bar{\mathbb{E}}^{1/p}
\biggl[ \frac{\|
W_0 \|_{\mathsf{X}, \infty}^p \| f_i
\|_{\mathsf{X}^{n+1}, \infty}^p }{\{\mu g_{}(Y_0)\}^p} \biggr] +
\bar{\mathbb{E}}^{1/p} [ \| f_i \|_{\mathsf
{X}^{n+1}, \infty}^p ] \Biggr\}.
\end{eqnarray*}
We finish the proof by substituting the bounds of assumption
(A4) into the expression above and
summing up.
The proof of the second part follows similarly.
\end{pf}

\begin{pf*}{Proof of Proposition \protect\ref{propositionrandobs}} The
proof of the first part follows by applying Proposition
\ref{propositionrandobsfi} and the bound (\ref{eqbackwardmixing})
to the decomposition
\begin{eqnarray*}
\| \widehat{\gamma}^{N,\Delta_n}_n - \gamma_n \|_p
&\leq&\sum_{k=0}^{n-1} \| \widehat{\phi}_{\nu, 0:k(\Delta
_n)|k(\Delta_n)}^N s_k -
\phi_{\nu, 0:k(\Delta_n)|k(\Delta_n)} s_k \|_p \\
&&{}+ \sum_{k=0}^{n - \Delta_n} \| \phi_{\nu, 0:k + \Delta_n|k +
\Delta_n} s_k
-\phi_{\nu, 0:n|n}s_k \|_p.
\end{eqnarray*}
The second part is proved in a similar manner.
\end{pf*}
\end{appendix}

\section*{Acknowledgements} This work was supported by a~grant from the
Swedish Foundation for Strategic Research, a~French government
overseas student grant and a grant from the French National Agency
for Research (ANR-2005 ADAP'MC project). The authors are grateful to
the anonymous referees who provided sensible comments on our results
that improved the presentation of the paper.

\printhistory

\end{document}